\documentclass[11pt,reqno]{amsart}

\usepackage[T1]{fontenc}
\usepackage{palatino}
\usepackage{amsmath,amssymb,amsthm,mathtools}
\usepackage{amscd}
\usepackage[margin=1.1in]{geometry}
\usepackage{float}
\usepackage{hyperref}
\usepackage[nameinlink,capitalize]{cleveref}
\usepackage{placeins}
\usepackage{tikz}
\usetikzlibrary{arrows.meta,intersections}
\usetikzlibrary{decorations.markings}
\usetikzlibrary{calc}

\hypersetup{colorlinks=true,linkcolor=blue,citecolor=blue,urlcolor=blue}

\theoremstyle{plain}
\newtheorem{theorem}{Theorem}[section]
\newtheorem{proposition}[theorem]{Proposition}
\newtheorem{corollary}[theorem]{Corollary}
\newtheorem{lemma}[theorem]{Lemma}

\newtheorem{conjecture}[theorem]{Conjecture}

\theoremstyle{definition}
\newtheorem{remark}[theorem]{Remark}

\newcommand{\Z}{\mathbb Z}

\title[Legendrian surgery on $ST^*K$ and infinitely many Stein fillings]{Legendrian surgery on the unit cotangent bundle of the Klein bottle and infinitely many Stein fillings}

\author{Burak Ozbagci}

\begin{document}
\begin{abstract}
We show that a single Legendrian surgery on the canonical contact structure of
the unit cotangent bundle $ST^*K$ of the Klein bottle $K$ produces a planar
contact $3$--manifold $(Y,\xi)$ with infinitely many pairwise non-homeomorphic
Stein fillings.  The fillings are obtained from the disk cotangent bundle
$DT^*K$ by attaching a single Weinstein $2$--handle along the images of one
Legendrian knot under the powers of a contactomorphism of
$(ST^*K,\xi_{\mathrm{can}})$, no positive power of which extends over $DT^*K$.
The construction takes place in a planar open book of
$(ST^*K,\xi_{\mathrm{can}})$, and $(Y,\xi)$ is
contactomorphic to Baykur's recent example of a planar contact $3$--manifold
with infinitely many Stein fillings.  Since $(ST^*K,\xi_{\mathrm{can}})$ has a
unique minimal symplectic filling up to homeomorphism, finiteness of fillings is not preserved
under Legendrian surgery, even for planar contact $3$--manifolds.  The same
construction gives the planar contact $3$--manifolds with infinitely many
rational homology ball Stein fillings found earlier by the author.  Finally, we use Luttinger twists along the Klein bottle to construct contact
$3$--manifolds with infinitely many rational homology ball Stein fillings, without the use of
open books or Lefschetz fibrations.
\end{abstract}

\maketitle

\section{Introduction}

Baykur \cite{BaykurPlanar} constructed a planar contact $3$--manifold, supported
by the open book $(\Sigma_0^6,t_dB)$ with $B=t_ut_vt_wt_z$, where $u,v,w,z,d$ are
the curves of \cref{fig:baykur-curves}, and $t_c$ denotes the right-handed Dehn
twist along a simple closed curve $c$.  We use functional notation for products in
the mapping class group, so that $t_at_b$ means that we first apply $t_b$.  Here $\Sigma_0^6$ is the disk with holes
$\delta_1,\ldots,\delta_5$, numbered from left to right, and outer boundary $\delta_0$.  Baykur's Stein fillings $X^{\mathrm B}_n$,
$n\ge0$, are the Lefschetz fibrations of the positive factorizations
$t_d\,t_\alpha^nBt_\alpha^{-n}$, where $\alpha$ is the curve of
\cref{fig:oo-curves}, whose Dehn twist commutes with $B$, and they are pairwise non-homeomorphic since
$\pi_1(X^{\mathrm B}_n)\cong\Z/(2n+1)\rtimes_{-1}\Z$.  In
\cite{OzSmooth} we identified the smooth topology of
Baykur's example, showed that the Lefschetz fibration of the factorization $B$
is diffeomorphic to the disk cotangent bundle $DT^*K$ of the Klein bottle $K$,
and constructed further fillings and rational homology ball fillings of related
planar contact $3$--manifolds.

In this paper we show that Baykur's construction can be carried out directly on
the canonical contact structure of the unit cotangent bundle $ST^*K$, the
boundary of $DT^*K$, and that this explains the infinitely many fillings
geometrically.  Let $(\Sigma_0^6,B')$, with
$B'=t_{V_0}t_{V_1}t_{V_2}t_{V_3}$, be the planar open book of Oba and the author
\cite[Theorem~3.7]{ObaOzbagci} supporting $(ST^*K,\xi_{\mathrm{can}})$.  Its
Lefschetz fibration is Johns' exact symplectic Lefschetz fibration on $DT^*K$
\cite{Johns}, and its vanishing cycles $V_0, V_1, V_2,V_3$, in the notation of
\cite[Figure~1]{LiOzbagci}, are shown in \cref{fig:oo-curves}.  Let $\beta$ be the round curve enclosing $\delta_3$ and $\delta_4$, shown in
\cref{fig:oo-curves}, and let $\sigma_i$ denote the half twist interchanging
$\delta_i$ and $\delta_{i+1}$ (\cref{sec:openbook}).  We consider the curves
\[
d'=\sigma_2(\beta),\qquad \alpha'=\sigma_1\sigma_2(\beta),
\]
shown in \cref{fig:oo-curves}.   We use primes because $B'$, $d'$ and
$\alpha'$ correspond to Baykur's $B$, $d$ and $\alpha$: there is a diffeomorphism
$g$ of the page $\Sigma_0^6$ with $gBg^{-1}=B'$, $g(d)=d'$ and $g(\alpha)=\alpha'$
(\cref{sec:baykur}).

Our main result is the following.

\begin{theorem}\label[theorem]{thm:main}
\leavevmode
\begin{enumerate}
\item The Dehn twist $t_{\alpha'}$ commutes with $B'$, and it induces a
contactomorphism $\Psi$ of $(ST^*K,\xi_{\mathrm{can}})$.
\item Let $(Y,\xi)$ be the contact $3$--manifold obtained from
$(ST^*K,\xi_{\mathrm{can}})$ by Legendrian surgery along a page realization of
$d'$.  For every $n\ge0$, the Stein domain
\[
X_n=DT^*K\cup h_{\Psi^{-n}(d')},
\]
obtained by attaching a Weinstein $2$--handle along the image of this Legendrian
knot under $\Psi^{-n}$, is a Stein filling of $(Y,\xi)$, and
$\pi_1(X_n)\cong\Z/(2n+1)\rtimes_{-1}\Z$.  In particular the fillings $X_n$ are
pairwise non-homeomorphic, and no positive power of $\Psi$ extends over
$DT^*K$, even as a homeomorphism.
\item The contact $3$--manifold $(Y,\xi)$ is contactomorphic to Baykur's planar
contact $3$--manifold, and $X_n$ is diffeomorphic to Baykur's filling
$X^{\mathrm B}_n$ for every $n\ge0$.
\end{enumerate}
\end{theorem}

Here $DT^*K$ carries the Stein structure of the Lefschetz fibration of
$(\Sigma_0^6,B')$, transported by the diffeomorphism constructed in the proof,
which is compatible with the boundary contact structures.  We expect it to be
deformation equivalent to the canonical one (\cref{rem:GP}).  The $2$--handle is
attached with framing one less than the page framing.  In the
proof we show that $(Y,\xi)$ is supported by the open book $(\Sigma_0^6,t_{d'}B')$,
and that $X_n$ is the Lefschetz fibration of the positive factorization
$t_{d'}\,t_{\alpha'}^nB't_{\alpha'}^{-n}$ of $t_{d'}B'$.  In particular, $X_0$ is the
Lefschetz fibration of $t_{d'}B'=t_{d'}t_{V_0}t_{V_1}t_{V_2}t_{V_3}$.   The point of \cref{thm:main} is that it takes place on
$(ST^*K,\xi_{\mathrm{can}})$ itself.  Part~(1) follows from
$B'(\alpha')=\alpha'$, which we see by following $\alpha'$ through the four twists of
$B'$ (\cref{lem:commute}), and part~(2) uses Johns'
theorem and a computation of fundamental groups (\cref{lem:pi1}).  For
part~(3) we show that Baykur's factorization $B$ becomes $B'$ after three Hurwitz
moves and a diffeomorphism $g$ of the page with $g(d)=d'$ and $g(\alpha)=\alpha'$
(\cref{sec:baykur}).

The monodromy of the canonical open book of $ST^*K$ also satisfies a
lantern-type relation, which comes from an arrangement of five real lines with
two triple points (\cref{prop:lantern}):
\[
t_{V_0}t_{V_1}t_{V_2}t_{V_3}\,t_{\alpha'}t_{\beta'}
=t_{\delta_0}\,t_{\delta_1}^2t_{\delta_2}^2\,t_{\delta_3}\,t_{\delta_4}^2t_{\delta_5}^2,
\]
where $\beta'$ is the curve of \cref{fig:oo-curves}, disjoint from $\alpha'$.  It is
Baykur's generalized lantern relation \cite[Lemma~3]{BaykurPlanar} transported by
$g$ (\cref{rem:transport}).  Thus $B'$ is a product of boundary twists composed
with negative Dehn twists along the two disjoint curves $\alpha'$ and $\beta'$.
Since boundary twists are central, this gives another explanation of
\cref{thm:main}(1).

\smallskip
\noindent\emph{A geometric point of view on the fillings.}  Part~(2) replaces
the list of positive factorizations by a single geometric picture with three
ingredients: a fixed filling, a Legendrian knot, and a contact mapping class.
The fixed filling is $DT^*K$.  By Li and the author \cite[Corollary~1.2]{LiOzbagci}, every
minimal weak symplectic filling of $(ST^*K,\xi_{\mathrm{can}})$ is homeomorphic
to $DT^*K$.  The Legendrian knot is a page realization of $d'$, and the contact
mapping class is $\Psi$.  Since $\Psi$ is a contactomorphism, Legendrian surgery
along every translate produces the same contact manifold $(Y,\xi)$, and the
pairwise non-isomorphic fundamental groups show that the framed attaching knots
are inequivalent under boundary diffeomorphisms that extend over the fixed
filling.  From this point of view, the infinitely many fillings are produced by
the gap between the contact mapping class group of the boundary and the
subgroup of mapping classes that extend over the filling.  As a consequence, we obtain the following.

\begin{corollary}\label[corollary]{cor:finiteness}
There is a Stein fillable planar contact $3$--manifold with a unique minimal
weak symplectic filling, up to homeomorphism, such that a single Legendrian surgery
on it yields a contact $3$--manifold with infinitely many Stein fillings, up to
homeomorphism.  In particular, finiteness of minimal symplectic fillings is
not preserved under Legendrian surgery, even within planar contact
$3$--manifolds.
\end{corollary}

In fact $\Psi$ has
infinite order in the smooth mapping class group of $ST^*K$, and no nontrivial
power of $\Psi$ is isotopic to the restriction of a homeomorphism of $DT^*K$
(\cref{cor:mcg}).  

\smallskip
\noindent\emph{Rational homology ball fillings.}  The planar contact
$3$--manifolds $(Y_m,\xi_m)$ and $(Y'_m,\xi'_m)$, $m\ge1$, of
\cite[Theorem~1.1]{OzSmooth} are supported by the open books
$(\Sigma_0^6,t_{d_m}B)$ and $(\Sigma_0^6,t_{d'_m}B)$, which differ from Baykur's
only in the curve $d$.  (The prime in $d'_m$ and $Y'_m$ is the notation of
\cite{OzSmooth}, and it is unrelated to the primes in $B'$, $d'$ and $\alpha'$.)
The same construction applies to them.

\begin{corollary}\label[corollary]{cor:qhb}
For every $m\ge1$, the contact $3$--manifolds $(Y_m,\xi_m)$ and $(Y'_m,\xi'_m)$
are obtained from $(ST^*K,\xi_{\mathrm{can}})$ by a single Legendrian surgery,
along page realizations of $g(d_m)$ and $g(d'_m)$, where $g$ is the
diffeomorphism of the page of \cref{lem:curve-identities}, and their pairwise
non-homeomorphic rational homology ball Stein fillings $X_{m,n}$ and $X'_{m,n}$,
$n\ge0$, of \cite[Theorem~1.1]{OzSmooth} are diffeomorphic to
$DT^*K\cup h_{\Psi^{-n}(g(d_m))}$ and $DT^*K\cup h_{\Psi^{-n}(g(d'_m))}$,
respectively.  In particular, a single Legendrian surgery on a contact
$3$--manifold with a unique minimal weak symplectic filling, up to homeomorphism,
can produce a contact $3$--manifold with infinitely many pairwise
non-homeomorphic rational homology ball Stein fillings.
\end{corollary}

\smallskip
\noindent\emph{Rational homology ball fillings from Luttinger twists.}  The phenomenon of \cref{cor:qhb}
can also be obtained directly on $T^*K$, without reference to an open book or a Lefschetz
fibration.  Luttinger surgery along a Lagrangian Klein bottle \cite{N02,N09b,N09} is a four-dimensional
cut-and-paste operation: a Weinstein neighborhood, which is a copy of $DT^*K$, is cut out and
reglued along its boundary $ST^*K$.  The gluing map is a diffeomorphism $f_{0,k}$ of $ST^*K$, which
we call a Luttinger twist, following \cite{N09b}.  In \cref{sec:luttinger} we give explicit strict
contactomorphisms $f'_k$ of $(ST^*K,\xi_{\mathrm{can}})$ isotopic to the Luttinger twists $f_{0,k}$ and satisfying $f'_k=(f'_1)^k$ (\cref{lem:strict}).  Attaching a
Weinstein $2$--handle to $DT^*K$, with its canonical symplectic structure, along the iterates of a
Legendrian knot under the single contactomorphism $f'_1$ then gives infinitely many rational homology ball Stein fillings of the
same contact $3$--manifold.  In the following statement, $s$, $t$ and $c$ are the loops in $ST^*K$
fixed in \cref{subsec:twists}, where $t$ is a fiber of $ST^*K\to K$.  We write the subgroup generated by
$s$ and $t$ additively, and we use the presentation $\pi_1(K)=\langle r,b\mid brb^{-1}=r^{-1}\rangle$.

\begin{theorem}\label[theorem]{thm:intro-qhb}
Let $\Lambda$ be a Legendrian knot in $(ST^*K,\xi_{\mathrm{can}})$ whose free homotopy class is $(ps+qt)\,c^2$ for
some integers $p$ and $q$ with $q\neq0$, and let $(Y_\Lambda,\xi_\Lambda)$ be the result of Legendrian surgery
along $\Lambda$.  Then, for every $k\in\Z$, the Stein domain obtained by attaching a Weinstein $2$--handle
to $DT^*K$, with its canonical symplectic structure, along $(f'_1)^k(\Lambda)$ is a rational homology ball Stein filling of $(Y_\Lambda,\xi_\Lambda)$ with
fundamental group
\[
\langle r,b\mid brb^{-1}=r^{-1},\ r^{|p+kq|}b^2=1\rangle ,
\]
which is the dicyclic group of order $4|p+kq|$ if $p+kq\neq0$ and the infinite dihedral group if
$p+kq=0$.  In particular,
infinitely many of these fillings are pairwise non-homeomorphic.
\end{theorem}

We prove \cref{thm:intro-qhb} in \cref{sec:luttinger}, where it is restated as \cref{thm:qhb-canonical}.  Such Legendrian knots exist for every $p\in\Z$ and every nonzero $q\in\Z$
(\cref{rem:existence}),
and the attaching circles of the fillings in \cref{cor:qhb} have images of the same form in
$\pi_1(K)$ (\cref{rem:qhb-compare}).  The effect of Klein bottle Luttinger surgery on closed
symplectic $4$--manifolds is studied in \cite{AE,DHN,Smirnov}.  By \cite[Remark~4]{BaykurPlanar}, consecutive
fillings in Baykur's family are related by fibered Luttinger surgeries along a Lagrangian torus.
Under our identification of the Lefschetz fibrations, the corresponding torus is swept out by
$\alpha'$ over a loop enclosing the four critical values of $B'$.  The partial conjugation describing
this surgery uses the same page Dehn twist $t_{\alpha'}$ that induces $\Psi$.  Separately, the fillings
of \cref{thm:intro-qhb} are constructed using strict contactomorphisms $f'_k$ of $ST^*K$ that are
smoothly isotopic to the Luttinger twists $f_{0,k}$.

\smallskip
\noindent\emph{Outline.}  In \cref{sec:openbook} we recall the canonical contact structure
$\xi_{\mathrm{can}}$ on $ST^*K$ and its planar open book $(\Sigma_0^6,B')$, introduce the curves $d'$,
$\alpha'$ and $\beta'$, and prove that $t_{\alpha'}$ commutes with $B'$
(\cref{lem:commute}).  We also compute the fundamental groups needed later
(\cref{lem:pi1}).  In \cref{subsec:lantern} we state a lantern-type relation for $B'$
(\cref{prop:lantern}), which comes from an arrangement of real lines, and we show
that the second configuration of vanishing cycles left open in \cite{LiOzbagci}
occurs (\cref{rem:fillings-ST}).  In
\cref{sec:fillings} we prove parts~(1) and~(2) of \cref{thm:main} and
\cref{cor:finiteness}, extend the construction to arbitrary page curves
(\cref{prop:general}), and derive the consequence for the mapping class groups of
$ST^*K$ (\cref{cor:mcg}).  In \cref{sec:baykur} we relate Baykur's factorization
to $B'$ by Hurwitz moves and a diffeomorphism of the page, which proves
\cref{thm:main}(3), and we deduce \cref{cor:qhb}.  In
\cref{sec:diagrams} we give a Stein diagram for the underlying smooth manifold of a
related filling $W_0$, and we state a conjecture.  In \cref{sec:luttinger} we use the
Luttinger twists along the Klein bottle to obtain rational homology ball Stein fillings
without reference to an open book.  In \cref{app:proofs} we give two proofs of
\cref{prop:lantern}.

\medskip
\noindent\emph{Acknowledgments.}  The author thanks Stefan Nemirovski for asking whether
the contactomorphism $\Psi$ is related to the contactomorphisms used to perform Luttinger
surgery along Lagrangian Klein bottles, which led to \cref{sec:luttinger}.

\medskip
\noindent\emph{Declaration on the use of AI.}  In the preparation of this paper
the author used LLM tools for
drawing the figures in TikZ from the author's own diagrams, and for editing the
text to improve the presentation.  The author takes full responsibility for the
content of the paper.

\section{The open book of \texorpdfstring{$ST^*K$}{ST*K}}\label{sec:openbook}

\subsection{The canonical contact structure}\label{subsec:canonical}

Let $\Sigma$ be a closed surface,
orientable or not, and let $\pi\colon T^*\Sigma\to\Sigma$ be its cotangent bundle.
The canonical $1$--form $\lambda_{\mathrm{can}}$ on $T^*\Sigma$ is defined at a covector
$p\in T^*_q\Sigma$ by $\lambda_{\mathrm{can}}(v)=p(d\pi(v))$, and in local coordinates
$(q_1,q_2)$ on $\Sigma$, with the dual fiber coordinates $(p_1,p_2)$, it is
$\lambda_{\mathrm{can}}=p_1\,dq_1+p_2\,dq_2$.  Its differential
$\omega_{\mathrm{can}}=d\lambda_{\mathrm{can}}=dp_1\wedge dq_1+dp_2\wedge dq_2$ is a
symplectic form, so $\omega_{\mathrm{can}}\wedge\omega_{\mathrm{can}}$ is a volume form,
and $T^*\Sigma$ is canonically oriented even when $\Sigma$ is nonorientable.  Fix a Riemannian metric on $\Sigma$, and let
$DT^*\Sigma$ and $ST^*\Sigma$ be the bundles of covectors of norm at most one and
exactly one.  The Liouville vector field of $\lambda_{\mathrm{can}}$ is the fiberwise
radial vector field $p_1\partial_{p_1}+p_2\partial_{p_2}$, which is transverse to
$ST^*\Sigma$ and points outward.  Hence $\lambda_{\mathrm{can}}$ restricts to a contact
form on $ST^*\Sigma$, whose kernel is the canonical contact structure
$\xi_{\mathrm{can}}$, and $(DT^*\Sigma,\omega_{\mathrm{can}})$ is an exact symplectic filling
of $(ST^*\Sigma,\xi_{\mathrm{can}})$, with $ST^*\Sigma$ oriented as its boundary.  Up to
isotopy, $\xi_{\mathrm{can}}$ does not depend on the metric.  Moreover $DT^*\Sigma$
carries a Stein structure inducing $\xi_{\mathrm{can}}$ on its boundary.  For a nonorientable
surface $\Sigma$, an explicit Stein diagram of $DT^*\Sigma$ can be obtained from the
solution of \cite[Exercise~11.2.5(a)]{GS}: Figure~12.75 there gives Stein diagrams of
the disk bundles over $\mathbb{RP}^2$, and Figure~12.76 shows how to obtain the general
case by band summing.  This is worked out in \cite[Figure~2]{OzArch}, and for the
Klein bottle the resulting Stein diagram of $DT^*K$ is shown in \cref{fig:W0-stein},
after removing the green knot.

\subsection{A planar open book supporting the canonical contact structure}\label{subsec:openbook}

For the closed nonorientable surface $N_k$ of genus $k$, Oba and the author
\cite[Theorem~3.7]{ObaOzbagci} constructed a planar open book, with page the disk
with $2k+1$ holes and monodromy the product of the Dehn twists along $k+2$
curves $V_0,\ldots,V_{k+1}$, supporting the canonical contact structure on the
unit cotangent bundle $ST^*N_k$.  For the Klein bottle $K=N_2$ this is an open book
$(\Sigma_0^6,B')$ supporting $(ST^*K,\xi_{\mathrm{can}})$ (the monodromy $B'$ is
denoted $\phi_2$ there), with
\[
B'=t_{V_0}\,t_{V_1}\,t_{V_2}\,t_{V_3},
\]
see also \cite[Lemma~2.3]{LiOzbagci}.

The left of \cref{fig:lioz} shows the
curves of \cite[Figure~1]{LiOzbagci} for $k=2$, where we write $V_3$ for
$V_{k+1}$.  The holes labeled $0,1,2,3,4$ there are our
$\delta_3,\delta_2,\delta_1,\delta_4,\delta_5$, and the outer boundary,
labeled $5$, is $\delta_0$.  Fix a basepoint $\ast\in\delta_0$ and the free generators $x_1,\ldots,x_5$
of $\pi_1(\Sigma_0^6,\ast)$ shown in \cref{fig:oo-curves}, which are the
 same generators of \cite[Figure~2(h)]{BaykurPlanar}.  We read products of paths from left to right.
 Moving the holes as indicated on the right of
\cref{fig:lioz} places them in a row, and
carries the curves $V_0,\dots,V_3$ to those of \cref{fig:oo-curves}.  In terms of
the generators $x_1,\dots,x_5$, they are given by the words
\begin{equation}\label{eq:Vwords}
V_0=x_4x_3x_2,\qquad V_1=x_2x_1,\qquad V_2=x_5x_4,\qquad V_3=x_5x_4x_3x_1x_4^{-1},
\end{equation}
up to conjugation and inversion.

\begin{figure}[ht]
\centering
\begin{tikzpicture}[scale=0.27]
\begin{scope}[shift={(0,0)}]
\draw[thick] (0,0) circle (10.6);
\draw[thick] (0.00,5.30) circle (0.6);
\draw[thick] (-3.00,1.00) circle (0.6);
\draw[thick] (-3.00,-3.60) circle (0.6);
\draw[thick] (3.00,1.00) circle (0.6);
\draw[thick] (3.00,-3.60) circle (0.6);
\draw[violet,thick] (0,2.6) circle (4.4);
\draw[green!55!black,thick] (-3,-1.3) ellipse (1.35 and 3.3);
\draw[green!55!black,thick] (3,-1.3) ellipse (1.35 and 3.3);
\draw[orange!85!black,thick] plot coordinates{(0.00,7.90) (0.17,7.94) (0.36,7.98) (0.57,8.02) (0.80,8.07) (1.04,8.12) (1.29,8.16) (1.54,8.19) (1.80,8.22) (2.06,8.24) (2.31,8.24) (2.56,8.23) (2.80,8.20) (3.04,8.16) (3.29,8.11) (3.54,8.06) (3.79,7.99) (4.05,7.92) (4.30,7.83) (4.55,7.73) (4.79,7.62) (5.01,7.49) (5.23,7.34) (5.42,7.18) (5.60,7.00) (5.76,6.80) (5.91,6.57) (6.04,6.32) (6.16,6.06) (6.27,5.78) (6.37,5.49) (6.46,5.18) (6.54,4.87) (6.61,4.56) (6.68,4.24) (6.74,3.92) (6.80,3.60) (6.85,3.28) (6.89,2.94) (6.91,2.59) (6.93,2.23) (6.94,1.86) (6.95,1.49) (6.95,1.13) (6.94,0.76) (6.94,0.41) (6.92,0.06) (6.91,-0.28) (6.90,-0.60) (6.89,-0.91) (6.87,-1.22) (6.85,-1.53) (6.83,-1.83) (6.81,-2.12) (6.78,-2.41) (6.75,-2.69) (6.71,-2.96) (6.67,-3.22) (6.62,-3.46) (6.56,-3.69) (6.50,-3.90) (6.43,-4.10) (6.36,-4.28) (6.28,-4.45) (6.20,-4.61) (6.11,-4.75) (6.01,-4.89) (5.91,-5.01) (5.80,-5.13) (5.69,-5.23) (5.57,-5.33) (5.44,-5.42) (5.30,-5.50) (5.15,-5.57) (4.98,-5.64) (4.80,-5.70) (4.61,-5.74) (4.42,-5.78) (4.22,-5.81) (4.02,-5.83) (3.82,-5.84) (3.62,-5.85) (3.44,-5.84) (3.26,-5.82) (3.10,-5.80) (2.95,-5.76) (2.79,-5.72) (2.64,-5.66) (2.49,-5.59) (2.35,-5.51) (2.21,-5.42) (2.08,-5.33) (1.96,-5.23) (1.85,-5.13) (1.76,-5.02) (1.67,-4.91) (1.60,-4.80) (1.54,-4.68) (1.50,-4.55) (1.47,-4.41) (1.45,-4.26) (1.44,-4.10) (1.44,-3.94) (1.45,-3.79) (1.46,-3.63) (1.49,-3.49) (1.52,-3.35) (1.56,-3.22) (1.60,-3.10) (1.65,-2.99) (1.71,-2.90) (1.79,-2.80) (1.87,-2.71) (1.96,-2.63) (2.05,-2.56) (2.15,-2.48) (2.26,-2.42) (2.36,-2.36) (2.48,-2.30) (2.59,-2.25) (2.70,-2.20) (2.82,-2.16) (2.94,-2.14) (3.08,-2.12) (3.21,-2.11) (3.36,-2.11) (3.50,-2.11) (3.64,-2.10) (3.79,-2.10) (3.92,-2.09) (4.06,-2.07) (4.18,-2.04) (4.30,-2.00) (4.41,-1.95) (4.52,-1.90) (4.63,-1.85) (4.73,-1.80) (4.83,-1.74) (4.93,-1.68) (5.02,-1.60) (5.11,-1.51) (5.19,-1.41) (5.27,-1.29) (5.34,-1.16) (5.40,-1.00) (5.45,-0.82) (5.50,-0.61) (5.54,-0.39) (5.58,-0.14) (5.61,0.11) (5.63,0.38) (5.65,0.66) (5.67,0.93) (5.68,1.21) (5.69,1.48) (5.70,1.75) (5.70,2.00) (5.71,2.25) (5.71,2.51) (5.72,2.76) (5.73,3.03) (5.73,3.29) (5.72,3.54) (5.71,3.80) (5.70,4.04) (5.67,4.28) (5.63,4.50) (5.57,4.71) (5.50,4.90) (5.41,5.08) (5.31,5.25) (5.19,5.41) (5.05,5.56) (4.91,5.71) (4.76,5.84) (4.60,5.96) (4.44,6.07) (4.27,6.17) (4.11,6.26) (3.95,6.34) (3.80,6.40) (3.65,6.45) (3.49,6.50) (3.32,6.53) (3.15,6.56) (2.98,6.57) (2.81,6.57) (2.65,6.56) (2.48,6.53) (2.32,6.50) (2.17,6.45) (2.03,6.38) (1.90,6.30) (1.78,6.21) (1.67,6.10) (1.56,5.98) (1.46,5.84) (1.37,5.69) (1.28,5.53) (1.20,5.36) (1.11,5.17) (1.04,4.97) (0.96,4.76) (0.88,4.54) (0.80,4.30) (0.72,4.05) (0.64,3.77) (0.56,3.47) (0.49,3.16) (0.41,2.83) (0.34,2.49) (0.27,2.15) (0.20,1.81) (0.14,1.47) (0.09,1.14) (0.04,0.81) (0.00,0.50) (-0.03,0.19) (-0.05,-0.12) (-0.06,-0.43) (-0.07,-0.75) (-0.07,-1.06) (-0.07,-1.37) (-0.07,-1.67) (-0.06,-1.96) (-0.06,-2.24) (-0.07,-2.51) (-0.08,-2.77) (-0.10,-3.00) (-0.12,-3.22) (-0.15,-3.42) (-0.18,-3.61) (-0.20,-3.79) (-0.23,-3.96) (-0.27,-4.12) (-0.31,-4.27) (-0.35,-4.41) (-0.40,-4.54) (-0.46,-4.67) (-0.53,-4.79) (-0.60,-4.90) (-0.68,-5.01) (-0.78,-5.11) (-0.88,-5.20) (-0.99,-5.28) (-1.10,-5.36) (-1.22,-5.43) (-1.34,-5.49) (-1.47,-5.54) (-1.60,-5.59) (-1.73,-5.63) (-1.87,-5.67) (-2.00,-5.70) (-2.14,-5.73) (-2.28,-5.75) (-2.44,-5.76) (-2.60,-5.77) (-2.76,-5.78) (-2.92,-5.78) (-3.08,-5.77) (-3.24,-5.75) (-3.39,-5.72) (-3.54,-5.69) (-3.67,-5.65) (-3.80,-5.60) (-3.92,-5.54) (-4.04,-5.47) (-4.16,-5.39) (-4.27,-5.30) (-4.38,-5.20) (-4.48,-5.09) (-4.58,-4.98) (-4.66,-4.87) (-4.74,-4.75) (-4.81,-4.64) (-4.86,-4.52) (-4.90,-4.40) (-4.93,-4.28) (-4.94,-4.15) (-4.94,-4.01) (-4.93,-3.86) (-4.91,-3.72) (-4.89,-3.57) (-4.85,-3.42) (-4.81,-3.28) (-4.76,-3.15) (-4.71,-3.02) (-4.66,-2.90) (-4.60,-2.80) (-4.54,-2.71) (-4.47,-2.62) (-4.39,-2.54) (-4.30,-2.47) (-4.21,-2.40) (-4.12,-2.34) (-4.02,-2.29) (-3.92,-2.24) (-3.81,-2.20) (-3.71,-2.16) (-3.60,-2.13) (-3.50,-2.10) (-3.39,-2.08) (-3.28,-2.08) (-3.16,-2.09) (-3.04,-2.10) (-2.91,-2.12) (-2.79,-2.14) (-2.66,-2.17) (-2.54,-2.19) (-2.42,-2.21) (-2.31,-2.21) (-2.20,-2.21) (-2.10,-2.20) (-2.01,-2.18) (-1.91,-2.16) (-1.82,-2.15) (-1.74,-2.13) (-1.65,-2.10) (-1.58,-2.07) (-1.50,-2.03) (-1.43,-1.97) (-1.36,-1.91) (-1.30,-1.82) (-1.25,-1.72) (-1.20,-1.60) (-1.16,-1.45) (-1.12,-1.28) (-1.09,-1.10) (-1.07,-0.89) (-1.05,-0.67) (-1.03,-0.44) (-1.02,-0.20) (-1.01,0.04) (-1.01,0.29) (-1.00,0.53) (-1.00,0.77) (-1.00,1.00) (-1.00,1.23) (-1.01,1.46) (-1.03,1.70) (-1.04,1.93) (-1.07,2.17) (-1.09,2.41) (-1.11,2.65) (-1.13,2.89) (-1.15,3.12) (-1.17,3.35) (-1.19,3.58) (-1.20,3.80) (-1.21,4.02) (-1.21,4.24) (-1.22,4.46) (-1.22,4.67) (-1.22,4.89) (-1.23,5.10) (-1.22,5.31) (-1.22,5.50) (-1.22,5.69) (-1.21,5.87) (-1.21,6.04) (-1.20,6.20) (-1.19,6.34) (-1.19,6.48) (-1.18,6.60) (-1.18,6.71) (-1.17,6.82) (-1.16,6.92) (-1.15,7.01) (-1.13,7.10) (-1.11,7.18) (-1.08,7.25) (-1.05,7.33) (-1.00,7.40) (-0.95,7.47) (-0.90,7.53) (-0.85,7.58) (-0.79,7.63) (-0.73,7.67) (-0.66,7.71) (-0.58,7.74) (-0.50,7.77) (-0.40,7.81) (-0.28,7.84) (-0.15,7.87)} -- cycle;
\node[violet] at (-4.2,6.0) {\small$V_0$};
\node[green!55!black] at (-5.2,-1.3) {\small$V_1$};
\node[green!55!black] at (1.75,2.35) {\small$V_2$};
\node[orange!85!black] at (7.90,1.49) {\small$V_3$};
\end{scope}
\begin{scope}[shift={(24,0)}]
\node at (8.3,-8.6) {\small$\delta_0$};
\draw[thick] (0,0) circle (10.6);
\draw[thick] (0.00,5.30) circle (0.6);
\node at (1.50,5.30) {\scriptsize$\delta_{3}$};
\draw[thick] (-3.00,1.00) circle (0.6);
\node at (-3.00,-0.40) {\scriptsize$\delta_{2}$};
\draw[thick] (-3.00,-3.60) circle (0.6);
\node at (-3.00,-5.00) {\scriptsize$\delta_{1}$};
\draw[thick] (3.00,1.00) circle (0.6);
\node at (3.00,-0.40) {\scriptsize$\delta_{4}$};
\draw[thick] (3.00,-3.60) circle (0.6);
\node at (3.00,-5.00) {\scriptsize$\delta_{5}$};
\draw[blue!70!black,thick,->] (0,4.5) -- (0,2.0);
\draw[blue!70!black,thick,->] (-3.5,-4.3) .. controls (-6.5,-6) and (-8.6,-2) .. (-8.3,1.0);
\draw[blue!70!black,thick,->] (3.5,-4.3) .. controls (6.5,-6) and (8.6,-2) .. (8.3,1.0);
\end{scope}
\end{tikzpicture}
\caption{Left: the curves $V_0,\dots,V_3$ of \cite[Figure~1]{LiOzbagci}
for $k=2$, with $V_3=V_{k+1}$.  Right: moving the holes along the arrows
places them in a row.}
\label{fig:lioz}
\end{figure}
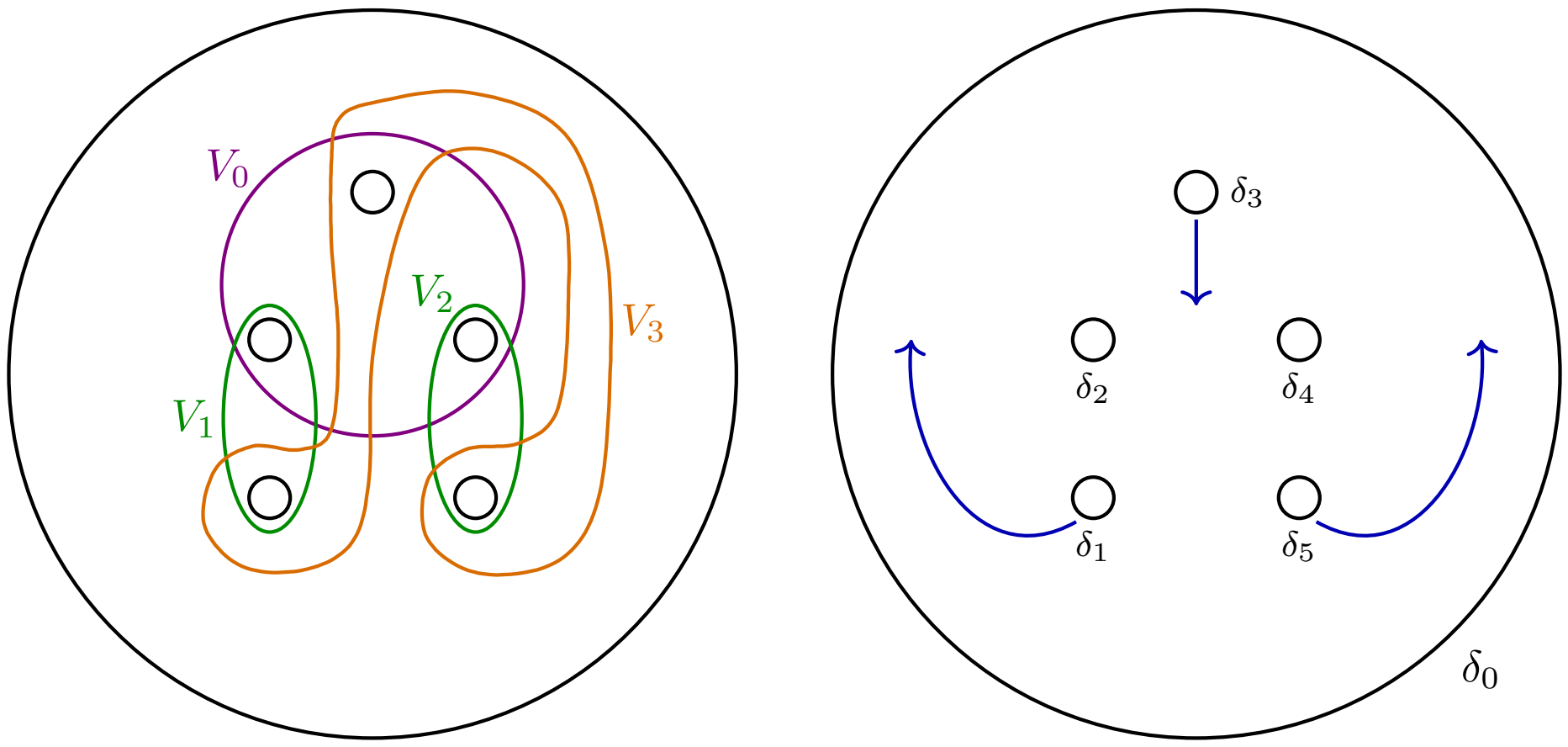

\begin{figure}[ht]
\centering
\makebox[\textwidth][c]{%
\begin{tikzpicture}[scale=0.3]
\begin{scope}[shift={(0.0,0.0)}]
\draw[thick,rounded corners=10pt] (0,-2.75) rectangle (15.600000000000001,2.75);
\draw[thick] (2.60,0.00) circle (0.35);
\draw[thick] (5.20,0.00) circle (0.35);
\draw[thick] (7.80,0.00) circle (0.35);
\draw[thick] (10.40,0.00) circle (0.35);
\draw[thick] (13.00,0.00) circle (0.35);
\draw[violet,thick] plot[smooth cycle,tension=0.7] coordinates{(5.14,-0.69) (5.01,-0.67) (4.90,-0.63) (4.80,-0.57) (4.72,-0.50) (4.65,-0.42) (4.59,-0.33) (4.54,-0.23) (4.51,-0.10) (4.50,0.03) (4.52,0.16) (4.56,0.28) (4.61,0.38) (4.68,0.46) (4.76,0.54) (4.85,0.60) (4.95,0.65) (5.07,0.69) (5.21,0.70) (5.37,0.70) (5.52,0.70) (5.68,0.70) (5.83,0.70) (5.99,0.70) (6.14,0.70) (6.30,0.70) (6.46,0.70) (6.61,0.70) (6.77,0.70) (6.92,0.70) (7.08,0.70) (7.23,0.70) (7.39,0.70) (7.55,0.70) (7.70,0.70) (7.86,0.70) (8.01,0.70) (8.17,0.70) (8.32,0.70) (8.48,0.70) (8.64,0.70) (8.79,0.70) (8.95,0.70) (9.10,0.70) (9.26,0.70) (9.41,0.70) (9.57,0.70) (9.73,0.70) (9.88,0.70) (10.04,0.70) (10.19,0.70) (10.35,0.70) (10.49,0.69) (10.62,0.66) (10.73,0.62) (10.82,0.56) (10.90,0.48) (10.97,0.40) (11.03,0.30) (11.07,0.19) (11.09,0.07) (11.09,-0.07) (11.07,-0.19) (11.03,-0.30) (10.97,-0.40) (10.90,-0.48) (10.82,-0.56) (10.73,-0.62) (10.62,-0.66) (10.49,-0.69) (10.35,-0.70) (10.19,-0.70) (10.04,-0.70) (9.88,-0.70) (9.73,-0.70) (9.57,-0.70) (9.41,-0.70) (9.26,-0.70) (9.10,-0.70) (8.95,-0.70) (8.79,-0.70) (8.64,-0.70) (8.48,-0.70) (8.32,-0.70) (8.17,-0.70) (8.01,-0.70) (7.86,-0.70) (7.70,-0.70) (7.55,-0.70) (7.39,-0.70) (7.23,-0.70) (7.08,-0.70) (6.92,-0.70) (6.77,-0.70) (6.61,-0.70) (6.46,-0.70) (6.30,-0.70) (6.14,-0.70) (5.99,-0.70) (5.83,-0.70) (5.68,-0.70) (5.52,-0.70) (5.37,-0.70) (5.21,-0.70)};
\node[violet] at (3.85,1.40) {$V_0$};
\end{scope}
\begin{scope}[shift={(16.8,0.0)}]
\draw[thick,rounded corners=10pt] (0,-2.75) rectangle (15.600000000000001,2.75);
\draw[thick] (2.60,0.00) circle (0.35);
\draw[thick] (5.20,0.00) circle (0.35);
\draw[thick] (7.80,0.00) circle (0.35);
\draw[thick] (10.40,0.00) circle (0.35);
\draw[thick] (13.00,0.00) circle (0.35);
\draw[green!55!black,thick] plot[smooth cycle,tension=0.7] coordinates{(2.53,-0.69) (2.45,-0.68) (2.37,-0.66) (2.30,-0.63) (2.23,-0.59) (2.17,-0.55) (2.11,-0.50) (2.07,-0.45) (2.02,-0.39) (1.98,-0.33) (1.95,-0.26) (1.93,-0.18) (1.91,-0.10) (1.90,-0.00) (1.91,0.09) (1.92,0.17) (1.95,0.25) (1.98,0.32) (2.02,0.39) (2.06,0.44) (2.11,0.50) (2.16,0.54) (2.22,0.59) (2.29,0.62) (2.36,0.66) (2.44,0.68) (2.53,0.69) (2.63,0.70) (2.73,0.70) (2.83,0.70) (2.94,0.70) (3.04,0.70) (3.14,0.70) (3.25,0.70) (3.35,0.70) (3.46,0.70) (3.56,0.70) (3.66,0.70) (3.77,0.70) (3.87,0.70) (3.97,0.70) (4.08,0.70) (4.18,0.70) (4.29,0.70) (4.39,0.70) (4.49,0.70) (4.60,0.70) (4.70,0.70) (4.80,0.70) (4.91,0.70) (5.01,0.70) (5.12,0.70) (5.22,0.70) (5.31,0.69) (5.40,0.67) (5.47,0.64) (5.54,0.61) (5.60,0.57) (5.66,0.52) (5.71,0.47) (5.76,0.42) (5.80,0.36) (5.83,0.29) (5.86,0.22) (5.88,0.14) (5.89,0.05) (5.90,-0.04) (5.88,-0.13) (5.86,-0.21) (5.84,-0.28) (5.80,-0.35) (5.76,-0.41) (5.72,-0.47) (5.67,-0.52) (5.61,-0.57) (5.55,-0.61) (5.48,-0.64) (5.40,-0.67) (5.32,-0.69) (5.23,-0.70) (5.13,-0.70) (5.02,-0.70) (4.92,-0.70) (4.82,-0.70) (4.71,-0.70) (4.61,-0.70) (4.50,-0.70) (4.40,-0.70) (4.30,-0.70) (4.19,-0.70) (4.09,-0.70) (3.98,-0.70) (3.88,-0.70) (3.78,-0.70) (3.67,-0.70) (3.57,-0.70) (3.47,-0.70) (3.36,-0.70) (3.26,-0.70) (3.15,-0.70) (3.05,-0.70) (2.95,-0.70) (2.84,-0.70) (2.74,-0.70) (2.64,-0.70) (2.54,-0.69)};
\node[green!55!black] at (1.25,1.40) {$V_1$};
\end{scope}
\begin{scope}[shift={(33.6,0.0)}]
\draw[thick,rounded corners=10pt] (0,-2.75) rectangle (15.600000000000001,2.75);
\draw[thick] (2.60,0.00) circle (0.35);
\draw[thick] (5.20,0.00) circle (0.35);
\draw[thick] (7.80,0.00) circle (0.35);
\draw[thick] (10.40,0.00) circle (0.35);
\draw[thick] (13.00,0.00) circle (0.35);
\draw[green!55!black,thick] plot[smooth cycle,tension=0.7] coordinates{(10.33,-0.69) (10.25,-0.68) (10.17,-0.66) (10.09,-0.63) (10.03,-0.59) (9.97,-0.55) (9.91,-0.50) (9.86,-0.45) (9.82,-0.39) (9.78,-0.33) (9.75,-0.26) (9.73,-0.18) (9.71,-0.10) (9.70,-0.00) (9.71,0.09) (9.72,0.17) (9.75,0.25) (9.78,0.32) (9.82,0.38) (9.86,0.44) (9.91,0.50) (9.96,0.54) (10.02,0.59) (10.09,0.62) (10.16,0.65) (10.24,0.68) (10.33,0.69) (10.42,0.70) (10.53,0.70) (10.63,0.70) (10.73,0.70) (10.84,0.70) (10.94,0.70) (11.04,0.70) (11.15,0.70) (11.25,0.70) (11.36,0.70) (11.46,0.70) (11.56,0.70) (11.67,0.70) (11.77,0.70) (11.87,0.70) (11.98,0.70) (12.08,0.70) (12.19,0.70) (12.29,0.70) (12.39,0.70) (12.50,0.70) (12.60,0.70) (12.71,0.70) (12.81,0.70) (12.91,0.70) (13.01,0.70) (13.11,0.69) (13.19,0.67) (13.27,0.64) (13.34,0.61) (13.40,0.57) (13.46,0.53) (13.51,0.48) (13.56,0.42) (13.60,0.36) (13.63,0.29) (13.66,0.22) (13.68,0.14) (13.69,0.05) (13.70,-0.04) (13.68,-0.13) (13.66,-0.21) (13.64,-0.29) (13.60,-0.35) (13.56,-0.42) (13.52,-0.47) (13.46,-0.52) (13.41,-0.57) (13.34,-0.61) (13.28,-0.64) (13.20,-0.67) (13.12,-0.69) (13.02,-0.70) (12.92,-0.70) (12.82,-0.70) (12.72,-0.70) (12.61,-0.70) (12.51,-0.70) (12.40,-0.70) (12.30,-0.70) (12.20,-0.70) (12.09,-0.70) (11.99,-0.70) (11.89,-0.70) (11.78,-0.70) (11.68,-0.70) (11.57,-0.70) (11.47,-0.70) (11.37,-0.70) (11.26,-0.70) (11.16,-0.70) (11.05,-0.70) (10.95,-0.70) (10.85,-0.70) (10.74,-0.70) (10.64,-0.70) (10.54,-0.70) (10.43,-0.70) (10.34,-0.69)};
\node[green!55!black] at (14.35,1.40) {$V_2$};
\end{scope}
\begin{scope}[shift={(0.0,-6.3)}]
\draw[thick,rounded corners=10pt] (0,-2.75) rectangle (15.600000000000001,2.75);
\draw[thick] (2.60,0.00) circle (0.35);
\draw[thick] (5.20,0.00) circle (0.35);
\draw[thick] (7.80,0.00) circle (0.35);
\draw[thick] (10.40,0.00) circle (0.35);
\draw[thick] (13.00,0.00) circle (0.35);
\draw[orange!85!black,thick] plot[smooth cycle,tension=0.7] coordinates{(2.57,-2.35) (2.27,-2.27) (2.06,-2.09) (1.93,-1.84) (1.90,-1.48) (1.90,-1.10) (1.90,-0.72) (1.90,-0.34) (1.90,0.04) (1.99,0.34) (2.17,0.55) (2.43,0.68) (2.77,0.68) (3.03,0.55) (3.21,0.34) (3.30,0.04) (3.30,-0.34) (3.30,-0.72) (3.46,-0.95) (3.85,-0.95) (4.24,-0.95) (4.64,-0.95) (5.03,-0.95) (5.43,-0.95) (5.82,-0.95) (6.22,-0.95) (6.61,-0.95) (7.00,-0.95) (7.10,-0.67) (7.10,-0.29) (7.10,0.10) (7.10,0.48) (7.10,0.86) (7.10,1.24) (7.10,1.62) (7.17,1.94) (7.33,2.16) (7.57,2.31) (7.92,2.35) (8.31,2.35) (8.71,2.35) (9.10,2.35) (9.50,2.35) (9.89,2.35) (10.29,2.35) (10.68,2.35) (11.08,2.35) (11.47,2.35) (11.86,2.35) (12.26,2.35) (12.65,2.35) (13.04,2.35) (13.34,2.26) (13.55,2.08) (13.68,1.82) (13.70,1.46) (13.70,1.08) (13.70,0.70) (13.70,0.32) (13.69,-0.06) (13.60,-0.35) (13.42,-0.56) (13.15,-0.68) (12.81,-0.67) (12.55,-0.54) (12.38,-0.32) (12.30,-0.02) (12.30,0.36) (12.30,0.74) (12.12,0.95) (11.73,0.95) (11.33,0.95) (10.94,0.95) (10.55,0.95) (10.15,0.95) (9.76,0.95) (9.36,0.95) (8.97,0.95) (8.58,0.95) (8.50,0.65) (8.50,0.27) (8.50,-0.12) (8.50,-0.50) (8.50,-0.88) (8.50,-1.26) (8.50,-1.64) (8.43,-1.95) (8.26,-2.17) (8.01,-2.31) (7.66,-2.35) (7.27,-2.35) (6.87,-2.35) (6.48,-2.35) (6.08,-2.35) (5.69,-2.35) (5.29,-2.35) (4.90,-2.35) (4.50,-2.35) (4.11,-2.35) (3.71,-2.35) (3.32,-2.35) (2.93,-2.35)};
\node[orange!85!black] at (1.50,1.20) {$V_3$};
\end{scope}
\begin{scope}[shift={(16.8,-6.3)}]
\draw[thick,rounded corners=10pt] (0,-2.75) rectangle (15.600000000000001,2.75);
\draw[thick] (2.60,0.00) circle (0.35);
\draw[thick] (5.20,0.00) circle (0.35);
\draw[thick] (7.80,0.00) circle (0.35);
\draw[thick] (10.40,0.00) circle (0.35);
\draw[thick] (13.00,0.00) circle (0.35);
\draw[red!75!black,thick] plot[smooth cycle,tension=0.7] coordinates{(2.57,-2.35) (2.35,-2.30) (2.17,-2.20) (2.04,-2.06) (1.94,-1.88) (1.90,-1.65) (1.90,-1.38) (1.90,-1.11) (1.90,-0.84) (1.90,-0.57) (1.90,-0.30) (1.90,-0.03) (1.94,0.21) (2.03,0.40) (2.16,0.54) (2.33,0.64) (2.56,0.70) (2.80,0.67) (2.99,0.58) (3.13,0.45) (3.24,0.28) (3.29,0.07) (3.30,-0.20) (3.30,-0.47) (3.30,-0.74) (3.37,-0.94) (3.64,-0.95) (3.92,-0.95) (4.20,-0.95) (4.48,-0.95) (4.76,-0.95) (5.04,-0.95) (5.32,-0.95) (5.60,-0.95) (5.88,-0.95) (6.16,-0.95) (6.45,-0.95) (6.73,-0.95) (7.01,-0.95) (7.29,-0.95) (7.57,-0.95) (7.85,-0.95) (8.13,-0.95) (8.41,-0.95) (8.69,-0.95) (8.97,-0.95) (9.25,-0.95) (9.53,-0.95) (9.70,-0.85) (9.70,-0.58) (9.70,-0.31) (9.70,-0.04) (9.73,0.20) (9.82,0.39) (9.95,0.53) (10.12,0.64) (10.34,0.69) (10.59,0.67) (10.78,0.59) (10.93,0.46) (11.03,0.29) (11.09,0.08) (11.10,-0.19) (11.10,-0.46) (11.10,-0.73) (11.10,-1.00) (11.10,-1.27) (11.10,-1.54) (11.08,-1.79) (11.01,-1.99) (10.89,-2.15) (10.73,-2.26) (10.53,-2.34) (10.27,-2.35) (9.99,-2.35) (9.71,-2.35) (9.42,-2.35) (9.14,-2.35) (8.86,-2.35) (8.58,-2.35) (8.30,-2.35) (8.02,-2.35) (7.74,-2.35) (7.46,-2.35) (7.18,-2.35) (6.90,-2.35) (6.62,-2.35) (6.34,-2.35) (6.06,-2.35) (5.78,-2.35) (5.50,-2.35) (5.22,-2.35) (4.94,-2.35) (4.66,-2.35) (4.38,-2.35) (4.10,-2.35) (3.82,-2.35) (3.54,-2.35) (3.26,-2.35) (2.98,-2.35) (2.70,-2.35)};
\node[red!75!black] at (10.40,1.50) {$\alpha'$};
\end{scope}
\begin{scope}[shift={(33.6,-6.3)}]
\draw[thick,rounded corners=10pt] (0,-2.75) rectangle (15.600000000000001,2.75);
\draw[thick] (2.60,0.00) circle (0.35);
\draw[thick] (5.20,0.00) circle (0.35);
\draw[thick] (7.80,0.00) circle (0.35);
\draw[thick] (10.40,0.00) circle (0.35);
\draw[thick] (13.00,0.00) circle (0.35);
\draw[red!75!black,thick] plot[smooth cycle,tension=0.7] coordinates{(5.14,-0.69) (4.92,-0.64) (4.75,-0.53) (4.62,-0.39) (4.53,-0.20) (4.50,0.04) (4.50,0.31) (4.50,0.58) (4.50,0.85) (4.50,1.12) (4.50,1.39) (4.50,1.66) (4.54,1.89) (4.64,2.06) (4.77,2.20) (4.95,2.30) (5.18,2.35) (5.46,2.35) (5.74,2.35) (6.02,2.35) (6.30,2.35) (6.58,2.35) (6.86,2.35) (7.14,2.35) (7.42,2.35) (7.70,2.35) (7.98,2.35) (8.26,2.35) (8.54,2.35) (8.82,2.35) (9.10,2.35) (9.38,2.35) (9.66,2.35) (9.94,2.35) (10.22,2.35) (10.50,2.35) (10.78,2.35) (11.06,2.35) (11.35,2.35) (11.63,2.35) (11.91,2.35) (12.19,2.35) (12.47,2.35) (12.75,2.35) (13.02,2.35) (13.25,2.30) (13.43,2.20) (13.56,2.06) (13.66,1.88) (13.70,1.65) (13.70,1.38) (13.70,1.11) (13.70,0.84) (13.70,0.57) (13.70,0.30) (13.70,0.03) (13.66,-0.21) (13.57,-0.40) (13.44,-0.54) (13.27,-0.64) (13.04,-0.70) (12.80,-0.67) (12.61,-0.58) (12.47,-0.45) (12.36,-0.28) (12.31,-0.07) (12.30,0.20) (12.30,0.47) (12.30,0.74) (12.23,0.94) (11.96,0.95) (11.68,0.95) (11.40,0.95) (11.12,0.95) (10.84,0.95) (10.56,0.95) (10.28,0.95) (10.00,0.95) (9.72,0.95) (9.44,0.95) (9.15,0.95) (8.87,0.95) (8.59,0.95) (8.31,0.95) (8.03,0.95) (7.75,0.95) (7.47,0.95) (7.19,0.95) (6.91,0.95) (6.63,0.95) (6.35,0.95) (6.07,0.95) (5.90,0.85) (5.90,0.58) (5.90,0.31) (5.90,0.04) (5.87,-0.20) (5.78,-0.39) (5.65,-0.53) (5.48,-0.64) (5.26,-0.69)};
\node[red!75!black] at (5.20,-1.50) {$\beta'$};
\end{scope}
\begin{scope}[shift={(16.8,-12.6)}]
\draw[thick,rounded corners=10pt] (0,-2.75) rectangle (15.600000000000001,2.75);
\draw[thick] (2.60,0.00) circle (0.35);
\draw[thick] (5.20,0.00) circle (0.35);
\draw[thick] (7.80,0.00) circle (0.35);
\draw[thick] (10.40,0.00) circle (0.35);
\draw[thick] (13.00,0.00) circle (0.35);
\draw[blue!70!black,thick] plot[smooth cycle,tension=0.7] coordinates{(5.17,-2.35) (4.95,-2.30) (4.77,-2.20) (4.64,-2.06) (4.54,-1.88) (4.50,-1.65) (4.50,-1.38) (4.50,-1.11) (4.50,-0.84) (4.50,-0.57) (4.50,-0.30) (4.50,-0.03) (4.54,0.21) (4.63,0.40) (4.76,0.54) (4.93,0.64) (5.16,0.70) (5.40,0.67) (5.59,0.58) (5.73,0.45) (5.84,0.28) (5.89,0.07) (5.90,-0.20) (5.90,-0.47) (5.90,-0.74) (5.96,-0.94) (6.12,-0.95) (6.29,-0.95) (6.45,-0.95) (6.62,-0.95) (6.78,-0.95) (6.95,-0.95) (7.11,-0.95) (7.28,-0.95) (7.44,-0.95) (7.61,-0.95) (7.78,-0.95) (7.95,-0.95) (8.11,-0.95) (8.28,-0.95) (8.44,-0.95) (8.61,-0.95) (8.77,-0.95) (8.94,-0.95) (9.10,-0.95) (9.27,-0.95) (9.43,-0.95) (9.60,-0.95) (9.70,-0.85) (9.70,-0.58) (9.70,-0.31) (9.70,-0.04) (9.73,0.20) (9.82,0.39) (9.95,0.53) (10.12,0.64) (10.34,0.69) (10.59,0.67) (10.78,0.59) (10.93,0.46) (11.03,0.29) (11.09,0.08) (11.10,-0.19) (11.10,-0.46) (11.10,-0.73) (11.10,-1.00) (11.10,-1.27) (11.10,-1.54) (11.08,-1.79) (11.01,-1.99) (10.89,-2.15) (10.73,-2.26) (10.53,-2.34) (10.27,-2.35) (9.99,-2.35) (9.71,-2.35) (9.53,-2.35) (9.37,-2.35) (9.20,-2.35) (9.04,-2.35) (8.87,-2.35) (8.71,-2.35) (8.54,-2.35) (8.38,-2.35) (8.21,-2.35) (8.05,-2.35) (7.88,-2.35) (7.72,-2.35) (7.55,-2.35) (7.39,-2.35) (7.22,-2.35) (7.05,-2.35) (6.89,-2.35) (6.72,-2.35) (6.56,-2.35) (6.39,-2.35) (6.23,-2.35) (6.06,-2.35) (5.86,-2.35) (5.58,-2.35) (5.30,-2.35)};
\node[blue!70!black] at (5.20,1.55) {$d'$};
\end{scope}
\begin{scope}[shift={(0.0,-12.6)}]
\draw[thick,rounded corners=10pt] (0,-2.75) rectangle (15.600000000000001,2.75);
\draw[thick] (2.60,0.00) circle (0.35);
\draw[thick] (5.20,0.00) circle (0.35);
\draw[thick] (7.80,0.00) circle (0.35);
\draw[thick] (10.40,0.00) circle (0.35);
\draw[thick] (13.00,0.00) circle (0.35);
\draw[red!75!black,thick] plot[smooth cycle,tension=0.7] coordinates{(2.53,-0.69) (2.45,-0.68) (2.37,-0.66) (2.30,-0.63) (2.23,-0.59) (2.17,-0.55) (2.11,-0.50) (2.07,-0.45) (2.02,-0.39) (1.98,-0.33) (1.95,-0.26) (1.93,-0.18) (1.91,-0.10) (1.90,-0.00) (1.91,0.09) (1.92,0.17) (1.95,0.25) (1.98,0.32) (2.02,0.39) (2.06,0.44) (2.11,0.50) (2.16,0.54) (2.22,0.59) (2.29,0.62) (2.36,0.66) (2.44,0.68) (2.53,0.69) (2.63,0.70) (2.73,0.70) (2.83,0.70) (2.94,0.70) (3.04,0.70) (3.14,0.70) (3.25,0.70) (3.35,0.70) (3.46,0.70) (3.56,0.70) (3.66,0.70) (3.77,0.70) (3.87,0.70) (3.97,0.70) (4.08,0.70) (4.18,0.70) (4.29,0.70) (4.39,0.70) (4.49,0.70) (4.60,0.70) (4.70,0.70) (4.80,0.70) (4.91,0.70) (5.01,0.70) (5.12,0.70) (5.22,0.70) (5.31,0.69) (5.40,0.67) (5.47,0.64) (5.54,0.61) (5.60,0.57) (5.66,0.52) (5.71,0.47) (5.76,0.42) (5.80,0.36) (5.83,0.29) (5.86,0.22) (5.88,0.14) (5.89,0.05) (5.90,-0.04) (5.88,-0.13) (5.86,-0.21) (5.84,-0.28) (5.80,-0.35) (5.76,-0.41) (5.72,-0.47) (5.67,-0.52) (5.61,-0.57) (5.55,-0.61) (5.48,-0.64) (5.40,-0.67) (5.32,-0.69) (5.23,-0.70) (5.13,-0.70) (5.02,-0.70) (4.92,-0.70) (4.82,-0.70) (4.71,-0.70) (4.61,-0.70) (4.50,-0.70) (4.40,-0.70) (4.30,-0.70) (4.19,-0.70) (4.09,-0.70) (3.98,-0.70) (3.88,-0.70) (3.78,-0.70) (3.67,-0.70) (3.57,-0.70) (3.47,-0.70) (3.36,-0.70) (3.26,-0.70) (3.15,-0.70) (3.05,-0.70) (2.95,-0.70) (2.84,-0.70) (2.74,-0.70) (2.64,-0.70) (2.54,-0.69)};
\draw[red!75!black,thick] plot[smooth cycle,tension=0.7] coordinates{(7.73,-0.69) (7.65,-0.68) (7.57,-0.66) (7.50,-0.63) (7.43,-0.59) (7.37,-0.55) (7.32,-0.50) (7.27,-0.45) (7.22,-0.39) (7.18,-0.33) (7.15,-0.26) (7.13,-0.18) (7.11,-0.10) (7.10,-0.00) (7.11,0.09) (7.13,0.18) (7.15,0.25) (7.18,0.32) (7.22,0.39) (7.26,0.44) (7.31,0.50) (7.36,0.55) (7.42,0.59) (7.49,0.63) (7.56,0.66) (7.64,0.68) (7.73,0.69) (7.83,0.70) (7.93,0.70) (8.03,0.70) (8.14,0.70) (8.24,0.70) (8.35,0.70) (8.45,0.70) (8.55,0.70) (8.66,0.70) (8.76,0.70) (8.86,0.70) (8.97,0.70) (9.07,0.70) (9.18,0.70) (9.28,0.70) (9.38,0.70) (9.49,0.70) (9.59,0.70) (9.69,0.70) (9.80,0.70) (9.90,0.70) (10.01,0.70) (10.11,0.70) (10.21,0.70) (10.32,0.70) (10.42,0.70) (10.51,0.69) (10.60,0.67) (10.67,0.64) (10.74,0.61) (10.80,0.57) (10.86,0.52) (10.91,0.47) (10.96,0.42) (11.00,0.36) (11.03,0.29) (11.06,0.22) (11.08,0.14) (11.10,0.05) (11.10,-0.04) (11.09,-0.13) (11.07,-0.21) (11.04,-0.28) (11.00,-0.35) (10.96,-0.41) (10.92,-0.47) (10.87,-0.52) (10.81,-0.57) (10.75,-0.61) (10.68,-0.64) (10.60,-0.67) (10.52,-0.69) (10.43,-0.70) (10.33,-0.70) (10.22,-0.70) (10.12,-0.70) (10.02,-0.70) (9.91,-0.70) (9.81,-0.70) (9.71,-0.70) (9.60,-0.70) (9.50,-0.70) (9.39,-0.70) (9.29,-0.70) (9.19,-0.70) (9.08,-0.70) (8.98,-0.70) (8.87,-0.70) (8.77,-0.70) (8.67,-0.70) (8.56,-0.70) (8.46,-0.70) (8.36,-0.70) (8.25,-0.70) (8.15,-0.70) (8.04,-0.70) (7.94,-0.70) (7.84,-0.70) (7.74,-0.69)};
\node[red!75!black] at (3.90,1.25) {$\alpha$};
\node[red!75!black] at (9.10,1.25) {$\beta$};
\end{scope}
\begin{scope}[shift={(33.6,-12.6)}]
\draw[thick,rounded corners=10pt] (0,-2.75) rectangle (15.6,2.75);
\draw[thick] (2.60,0.00) circle (0.35);
\draw[blue!70!black] (7.80,-2.75) -- (3.351,-0.397);
\draw[blue!70!black,decoration={markings,mark=at position 0.327 with {\arrow{>}}},postaction={decorate}] (3.351,-0.397) arc (-27.87:332.13:0.85);
\node[blue!70!black] at (2.60,1.95) {\footnotesize$x_1$};
\draw[thick] (5.20,0.00) circle (0.35);
\draw[blue!70!black] (7.80,-2.75) -- (5.784,-0.618);
\draw[blue!70!black,decoration={markings,mark=at position 0.379 with {\arrow{>}}},postaction={decorate}] (5.784,-0.618) arc (-46.61:313.39:0.85);
\node[blue!70!black] at (5.20,1.95) {\footnotesize$x_2$};
\draw[thick] (7.80,0.00) circle (0.35);
\draw[blue!70!black] (7.80,-2.75) -- (7.800,-0.850);
\draw[blue!70!black,decoration={markings,mark=at position 0.500 with {\arrow{>}}},postaction={decorate}] (7.800,-0.850) arc (-90.00:270.00:0.85);
\node[blue!70!black] at (7.80,1.95) {\footnotesize$x_3$};
\draw[thick] (10.40,0.00) circle (0.35);
\draw[blue!70!black] (7.80,-2.75) -- (9.816,-0.618);
\draw[blue!70!black,decoration={markings,mark=at position 0.621 with {\arrow{>}}},postaction={decorate}] (9.816,-0.618) arc (-133.39:226.61:0.85);
\node[blue!70!black] at (10.40,1.95) {\footnotesize$x_4$};
\draw[thick] (13.00,0.00) circle (0.35);
\draw[blue!70!black] (7.80,-2.75) -- (12.249,-0.397);
\draw[blue!70!black,decoration={markings,mark=at position 0.673 with {\arrow{>}}},postaction={decorate}] (12.249,-0.397) arc (-152.13:207.87:0.85);
\node[blue!70!black] at (13.00,1.95) {\footnotesize$x_5$};
\fill (7.80,-2.75) circle (2.2pt); \node[below=-1pt] at (7.80,-2.75) {\footnotesize$\ast$};
\end{scope}
\end{tikzpicture}}
\caption{Top two rows: the vanishing cycles $V_0,\dots,V_3$ of $B'$, and the curves
$\alpha'$ and $\beta'$.  Bottom row: Baykur's curves $\alpha$ and $\beta$, the curve $d'$,
and the free generators $x_1,\dots,x_5$ of $\pi_1(\Sigma_0^6,\ast)$, with
$\ast\in\delta_0$.}
\label{fig:oo-curves}
\end{figure}
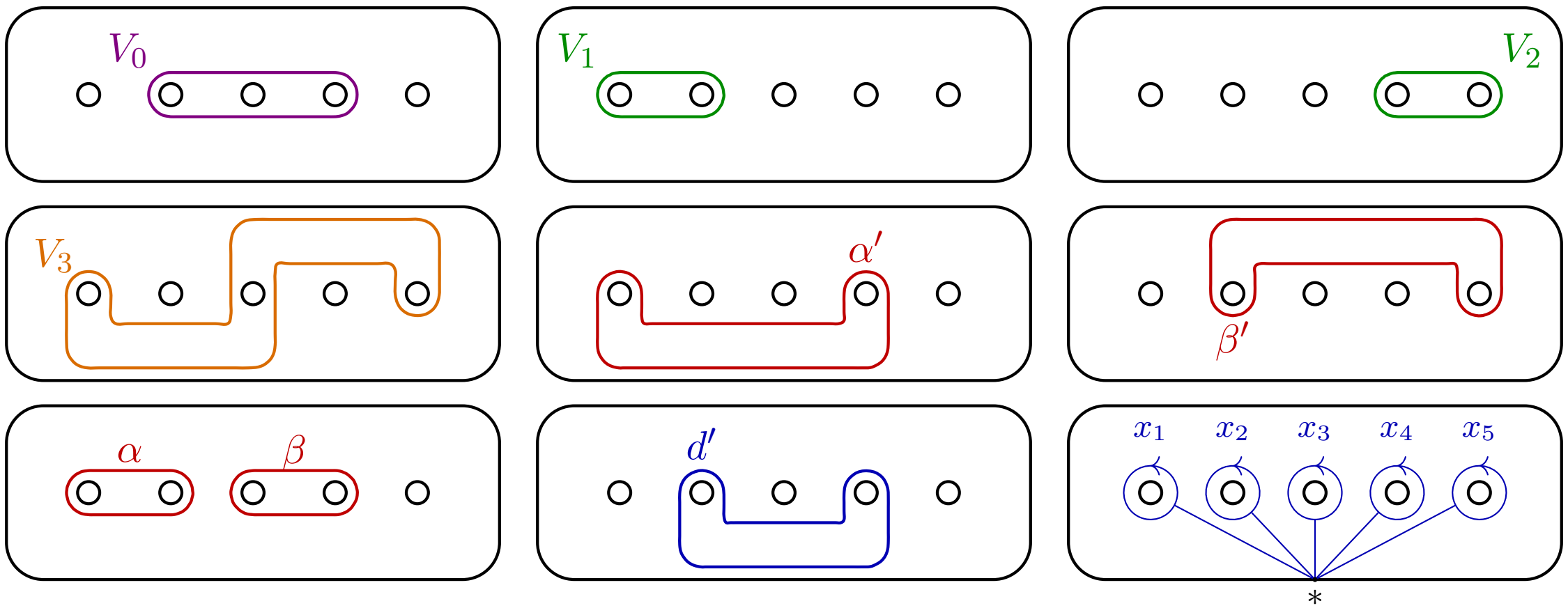

We also fix the half-twist convention on $\pi_1(\Sigma_0^6,\ast)$.  With our
convention for path products, the boundary of the two-hole disk containing $\delta_i$ and
$\delta_{i+1}$ is represented by $\gamma_i=x_{i+1}x_i$. For $1\le i\le4$, let $\sigma_i$ denote the clockwise half twist
interchanging $\delta_i$ and $\delta_{i+1}$, with the page oriented as in
\cref{fig:oo-curves}. Its action is
\begin{equation}\label{eq:artin}
\sigma_i\colon\quad x_i\mapsto x_i^{-1}x_{i+1}x_i,\qquad
x_{i+1}\mapsto x_i,\qquad x_j\mapsto x_j\ (j\ne i,i+1).
\end{equation}
Its inverse acts by
\begin{equation}\label{eq:artininv}
\sigma_i^{-1}\colon\quad x_i\mapsto x_{i+1},\qquad
x_{i+1}\mapsto x_{i+1}x_ix_{i+1}^{-1},\qquad
x_j\mapsto x_j\ (j\ne i,i+1).
\end{equation}
These formulas fix $\gamma_i$ and satisfy
$\sigma_i^2(x_j)=\gamma_i^{-1}x_j\gamma_i$ for $j=i,i+1$. Thus they are compatible
with the based boundary word, and $\sigma_i^2$ acts as the left-handed Dehn
twist along the boundary of the disk containing $\delta_i,\delta_{i+1}$.

A right-handed Dehn twist along the boundary of a disk containing the
consecutive holes $\delta_i,\ldots,\delta_j$ conjugates each of
$x_i,\ldots,x_j$ by the boundary word $x_j\cdots x_i$ and fixes the other
generators, since $\sigma_i^2$ is the left-handed twist.  In particular
$\beta=x_4x_3$, and \eqref{eq:artin} gives
\[
d'=\sigma_2(\beta)=x_4x_2,\qquad \alpha'=\sigma_1\sigma_2(\beta)=x_4x_1\sim x_1x_4 ,
\]
where $\sim$ denotes conjugacy.  We also use the curve $\beta'$ of
\cref{fig:oo-curves}, which is disjoint from $\alpha'$ and represents
$x_3^{-1}x_4^{-1}x_5x_4x_3x_2$ up to conjugacy.  Since freely homotopic essential simple closed
curves on a surface are isotopic \cite{Epstein}, identities between curves can be
checked on conjugacy classes in $\pi_1(\Sigma_0^6,\ast)$, up to inversion.

\begin{lemma}\label[lemma]{lem:commute}
$B'(\alpha')=\alpha'$ up to isotopy.  Hence $t_{\alpha'}$ commutes with $B'$.
\end{lemma}

\begin{figure}[ht]
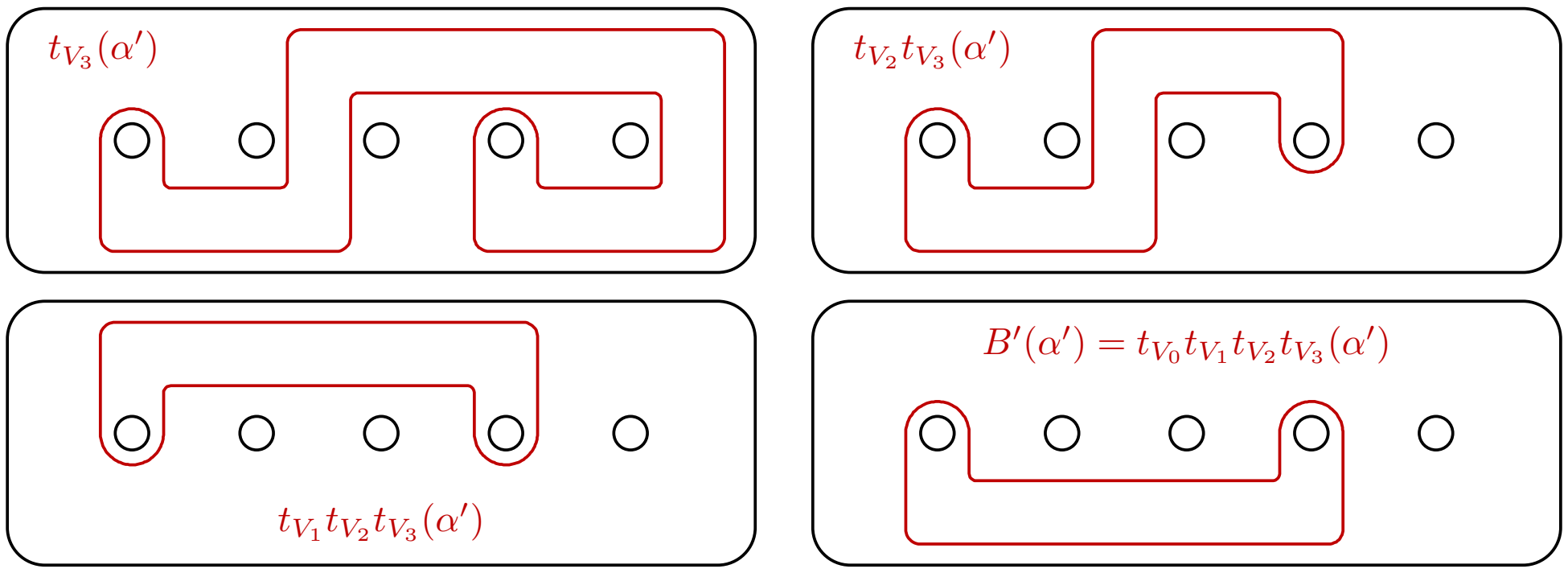

\centering
\resizebox{0.8\textwidth}{!}{
}
\caption{The images of $\alpha'$ under the four twists of $B'=t_{V_0}t_{V_1}t_{V_2}t_{V_3}$, applied from right to left (\cref{lem:commute}).  The last curve is $\alpha'$.}
\label{fig:Balpha}
\end{figure}

\begin{proof}
Since $B't_{\alpha'}B'^{-1}=t_{B'(\alpha')}$, the second statement follows from the
first.  Applying the Dehn twists $t_{V_3}$, $t_{V_2}$, $t_{V_1}$ and $t_{V_0}$ to
$\alpha'$ in turn gives the four curves of \cref{fig:Balpha}, each obtained from the
previous one by twisting along the corresponding curve of \cref{fig:oo-curves}.
The last of them is $\alpha'$.
\end{proof}

Let $X_{B'}$ denote the total space of the Lefschetz fibration with fiber
$\Sigma_0^6$ and vanishing cycles $V_0,V_1,V_2,V_3$.

\begin{lemma}\label[lemma]{lem:pi1}
The fundamental group of $X_{B'}$ is
$Q=\langle r,b\mid brb^{-1}=r^{-1}\rangle\cong\pi_1(K)$, where the generators of
$\pi_1(\Sigma_0^6,\ast)$ map to
\[
x_1\mapsto rb,\qquad x_2\mapsto b^{-1}r^{-1},\qquad x_3\mapsto r^{-1},\qquad
x_4\mapsto b,\qquad x_5\mapsto b^{-1}.
\]
For every $n\ge0$, the image of the curve $t_{\alpha'}^{-n}(d')$ in $Q$ is
conjugate to $r^{-(2n+1)}$.
\end{lemma}

\begin{proof}
The fundamental group of a Lefschetz fibration over the disk is the quotient
of the fundamental group of its fiber by the normal closure of its vanishing
cycles \cite[\S8.2]{GS}.  By \eqref{eq:Vwords}, the relation $V_1=1$ gives
$x_2=x_1^{-1}$, the relation $V_2=1$ gives $x_5=x_4^{-1}$, then $V_3=1$ gives
$x_3x_1x_4^{-1}=1$, i.e., $x_3=x_4x_1^{-1}$, and finally $V_0=1$ becomes
$x_4^2x_1^{-2}=1$.  So $\pi_1(X_{B'})=\langle x_1,x_4\mid x_1^2=x_4^2\rangle$.
With $r=x_1x_4^{-1}$ and $b=x_4$, so that $x_1=rb$, the relation $x_1^2=x_4^2$
becomes $rbrb=b^2$, which is equivalent to $brb^{-1}=r^{-1}$.  This gives the
presentation $Q$ and the images of the generators.  In $Q$ the element $b^2$ is
central, and $b^{-1}r^k=r^{-k}b^{-1}$ for every $k\in\Z$.

Let $\eta=\sigma_1\sigma_2$, so that $\alpha'=\eta(\beta)$ and
$t_{\alpha'}=\eta\,t_\beta\,\eta^{-1}$.  By \eqref{eq:artininv},
$\eta^{-1}(d')=\eta^{-1}(x_4x_2)=x_4x_3x_1x_3^{-1}$.  The twist $t_\beta^{-n}$
conjugates $x_3$ and $x_4$ by $\beta^{-n}$, where $\beta=x_4x_3$, so it maps
$x_4x_3x_1x_3^{-1}$ to $\beta x_1\beta^{-n}x_3^{-1}\beta^n$.  By \eqref{eq:artin},
$\eta(x_1)=x_1^{-1}x_2x_1$, $\eta(x_3)=x_1$ and $\eta(x_4)=x_4$, so $\eta(\beta)=x_4x_1$ and
\[
t_{\alpha'}^{-n}(d')=(x_4x_1)\,(x_1^{-1}x_2x_1)\,(x_4x_1)^{-n}\,x_1^{-1}\,(x_4x_1)^n .
\]
In $Q$ we have $x_4x_1\mapsto brb=r^{-1}b^2$, and $x_1^{-1}x_2x_1$ and $x_1^{-1}$
both map to $b^{-1}r^{-1}=rb^{-1}$.  Since $b^2$ is central, the image is
\[
(r^{-1}b^2)(rb^{-1})(r^nb^{-2n})(rb^{-1})(r^{-n}b^{2n})
=b^{-1}r^{n+1}b^{-1}r^{-n}b^2=r^{-(n+1)}b^{-2}r^{-n}b^2=r^{-(2n+1)} .\qedhere
\]
\end{proof}

\subsection{A lantern-type relation}\label{subsec:lantern}

In this subsection we state a lantern-type relation satisfied by the monodromy
$B'$.  It expresses $B'$ as a product of boundary twists composed with negative
Dehn twists along two disjoint curves $\alpha'$ and $\beta'$.  Since boundary twists
are central, this gives another reason why $t_{\alpha'}$ commutes with $B'$.  The
relation comes from an arrangement of five real lines with two triple points, by
Hironaka's construction of relations in the mapping class group from line
arrangements \cite{Hironaka}, and it is Baykur's generalized lantern relation
\cite[Lemma~3]{BaykurPlanar} transported by the diffeomorphism $g$ of
\cref{sec:baykur} (\cref{rem:transport}).

\begin{proposition}\label[proposition]{prop:lantern}
We have
\[
B'\,t_{\alpha'}t_{\beta'}=\Delta',\qquad
\Delta'=t_{\delta_0}\,t_{\delta_1}^2t_{\delta_2}^2\,t_{\delta_3}\,t_{\delta_4}^2t_{\delta_5}^2 .
\]
\end{proposition}

We give two proofs of \cref{prop:lantern} in \cref{app:proofs}.  The first follows
Baykur's proof of his generalized lantern relation via a line arrangement
\cite{BaykurPlanar,Hironaka}, and the second is a direct computation of the
action on the fundamental group of the page.

\begin{remark}\label[remark]{rem:fillings-ST}
By Wendl \cite{Wendl} and Niederkr\"uger--Wendl \cite{NW}, every minimal weak
symplectic filling of $(ST^*K,\xi_{\mathrm{can}})$ is symplectically deformation
equivalent to a Stein filling given by a positive allowable Lefschetz fibration,
that is, by a positive factorization of $B'$.  Li and the author showed
\cite[Lemma~2.8 and Remark~2.9]{LiOzbagci} that the vanishing
cycles of such a factorization enclose the same holes as $V_0,\ldots,V_3$, or
else the holes $\{\delta_1,\delta_2,\delta_3\}$, $\{\delta_3,\delta_4,\delta_5\}$,
$\{\delta_2,\delta_4\}$ and $\{\delta_1,\delta_5\}$, a possibility left open in
\cite[Remark~2.9(2)]{LiOzbagci}.  This second possibility does occur.  Let
$\sigma_{\beta'}$ be the half twist supported in the disk bounded by $\beta'$,
interchanging $\delta_2$ and $\delta_5$.  It fixes $\alpha'$ and $\beta'$, and it
exchanges $t_{\delta_2}$ and $t_{\delta_5}$, which have the same exponent in $\Delta'$,
so by \cref{prop:lantern} it commutes with $B'=\Delta't_{\beta'}^{-1}t_{\alpha'}^{-1}$.
Conjugating the factorization $t_{V_0}t_{V_1}t_{V_2}t_{V_3}$ by $\sigma_{\beta'}$ gives a
positive factorization of $B'$ whose vanishing cycles enclose the holes of the
second possibility.  Since it is a global conjugate, its Lefschetz fibration is
isomorphic to that of $t_{V_0}t_{V_1}t_{V_2}t_{V_3}$.
\end{remark}

\section{The infinitely many Stein fillings}\label{sec:fillings}

For a page curve $c$, i.e., a homologically nontrivial simple closed curve on a
page, a \emph{page realization} of $c$ is a Legendrian knot isotopic to $c$ on
a (possibly slightly perturbed) page, with contact framing equal to the page
framing.  Such a realization exists by Honda's Legendrian realization principle
\cite{Honda}.  Adding the twist $t_c$ to the monodromy of an open book
corresponds to Legendrian surgery along a page realization of $c$ \cite{AO}.

\begin{proof}[Proof of \cref{thm:main}(1) and (2)]
Oba and the author considered in \cite[\S3.2]{ObaOzbagci} Johns' exact
symplectic Lefschetz fibrations on the disk cotangent bundles of the
nonorientable surfaces $N_k$.  For the Klein bottle $K=N_2$, the one denoted
$\pi\colon(M,\omega)\to D^2$ there has fiber $\Sigma_0^6$ and ordered vanishing
cycles $(V_0,V_1,V_2,V_3)$, so $M=X_{B'}$, and its boundary carries the contact
structure $\xi_J=\ker(\lambda|_{\partial M})$, where $\lambda$ is the Liouville form of
$(M,\omega)$, supported by the induced open book $(\Sigma_0^6,B')$.  By
Johns' theorem \cite[Theorem~1.1]{Johns}, as recalled in \cite[\S3]{ObaOzbagci},
$(M,\omega)$ is, after smoothing its corners, conformally exact symplectomorphic
to $(DT^*K,\omega_{\mathrm{can}})$ by a map $\Phi$ with
$\Phi^*\lambda'=c\lambda+d\theta$, $c>0$, for a function $\theta$.  Here $\lambda'$ is
the Liouville form of $DT^*K$, so that both $\lambda$ and
$\Phi^*\lambda'$ restrict to contact forms on $\partial M$.  The forms
$\alpha_t=(c\lambda+t\,d\theta)|_{\partial M}$, $t\in[0,1]$, satisfy
$d\alpha_t=c\,d\lambda|_{\partial M}$, so $\alpha_t\wedge d\alpha_t$ is a convex
combination of the positive volume forms $\alpha_0\wedge d\alpha_0$ and
$\alpha_1\wedge d\alpha_1$, and each $\alpha_t$ is a contact form.  By Gray's
theorem, after an isotopy supported in a collar of the boundary, $\Phi$
restricts to a contactomorphism from $(\partial M,\xi_J)$ to
$(ST^*K,\xi_{\mathrm{can}})$.

The total space $X_{B'}$ also carries the Stein structure of the positive allowable
Lefschetz fibration with vanishing cycles $V_0,\ldots,V_3$ \cite{AO,LP}, and its
boundary contact structure is also supported by $(\Sigma_0^6,B')$.  By Giroux's
uniqueness theorem \cite{GirouxICM} the two boundary contact structures are
isotopic, and by Gray's theorem this isotopy is realized by an isotopy of
$\partial M$, which we extend over a collar.  Composing the resulting
diffeomorphism of $X_{B'}$, isotopic to the identity, with $\Phi$, we obtain a
diffeomorphism from $X_{B'}$ to $DT^*K$ which restricts to a contactomorphism from
the boundary of the Stein domain $X_{B'}$ to $(ST^*K,\xi_{\mathrm{can}})$.  We
transport the Stein structure of the Lefschetz fibration to $DT^*K$ by this
diffeomorphism.  From now on we write $\Phi$ for this diffeomorphism, and
$(Y_{B'},\xi_{B'})$ for the boundary of the Stein domain $X_{B'}$, so that
$\Phi|_{Y_{B'}}$ is a contactomorphism from $(Y_{B'},\xi_{B'})$ to
$(ST^*K,\xi_{\mathrm{can}})$.

We construct the contactomorphism $\widetilde\Psi$ of $(Y_{B'},\xi_{B'})$ induced by
$t_{\alpha'}$.  A mapping class is an isotopy class of diffeomorphisms, and to
build mapping tori we fix representatives.  We denote the representatives of
$B'$ and $t_{\alpha'}$ by the same letters, and we choose them equal to the
identity near the boundary, with $t_{\alpha'}$ supported in an annular
neighborhood of $\alpha'$.  The map $(x,s)\mapsto(t_{\alpha'}(x),s)$ is a
diffeomorphism from the mapping torus of $B'$ to the mapping torus of
$t_{\alpha'}B't_{\alpha'}^{-1}$.  Since $t_{\alpha'}B't_{\alpha'}^{-1}=B'$ as
mapping classes by \cref{lem:commute}, these two diffeomorphisms are isotopic rel
$\partial\Sigma_0^6$, and such an isotopy
induces a diffeomorphism between their mapping tori carrying pages to pages.
The composite, extended by the identity on the binding, is a diffeomorphism
$\widetilde\Psi_0$ of $Y_{B'}$ which preserves the open book, mapping each page to a page
and the binding to itself and preserving their orientations.  We choose the
isotopy so that the identification of the mapping tori is the identity on a
reference page, so that $\widetilde\Psi_0$ restricts to $t_{\alpha'}$ on that
page.  Hence $(\widetilde\Psi_0)_*\xi_{B'}$ is again supported by $(\Sigma_0^6,B')$, so it
is isotopic to $\xi_{B'}$ by Giroux's uniqueness theorem \cite{GirouxICM}.
Composing $\widetilde\Psi_0$ with the diffeomorphism given by Gray's theorem along such an
isotopy yields a contactomorphism $\widetilde\Psi$ of $(Y_{B'},\xi_{B'})$, smoothly
isotopic to $\widetilde\Psi_0$.  Let $\Psi$ be the contactomorphism of
$(ST^*K,\xi_{\mathrm{can}})$ obtained by transporting $\widetilde\Psi$ by $\Phi|_{Y_{B'}}$.
This proves part~(1).

For part~(2), let $L$ be a page realization of $d'$, and let $j$ denote the
framed embedding, into $Y_{B'}$, of a neighborhood of a copy of $d'$ on the
reference page, with framing one less than the page framing.  Since $\widetilde\Psi_0$ maps
pages to pages and preserves page framings, $\widetilde\Psi_0^{-n}\circ j$ is the framed
embedding of a copy of $t_{\alpha'}^{-n}(d')$ on a page, again with framing one
less than the page framing.  Since $\widetilde\Psi$ is smoothly isotopic to $\widetilde\Psi_0$, the
framed embedding $\widetilde\Psi^{-n}\circ j$ is framed isotopic to $\widetilde\Psi_0^{-n}\circ j$.
Attaching a $2$--handle to $X_{B'}$ along $\widetilde\Psi^{-n}\circ j$ therefore gives, as
a smooth manifold, the Lefschetz fibration with vanishing cycles
$t_{\alpha'}^{-n}(d'),V_0,\ldots,V_3$, which is the Lefschetz fibration of the
positive factorization
\[
P'_n=t_{d'}\,t_{\alpha'}^nB't_{\alpha'}^{-n}
\]
of $t_{d'}B'$, conjugated by $t_{\alpha'}^{-n}$.  Transporting by $\Phi$, and
writing again $L$ and $j$ for the transports of $L$ and $j$ to $ST^*K$, the
Stein domain $X_n=DT^*K\cup h_{\Psi^{-n}(L)}$ is obtained by attaching a
Weinstein $2$--handle along the Legendrian knot $\Psi^{-n}(L)$, whose contact
framing is the image of the page framing of $L$.  So the handle is attached
along $\Psi^{-n}\circ j$, and $X_n$ is diffeomorphic to the Lefschetz fibration
of $P'_n$.  Since $\Psi$ is a contactomorphism, Legendrian surgery along
$\Psi^{-n}(L)$ gives the same contact manifold as Legendrian surgery along $L$,
namely $(Y,\xi)$, which is supported by $(\Sigma_0^6,t_{d'}B')$.  Hence each
$X_n$ is a Stein filling of $(Y,\xi)$.

By \cref{lem:pi1}, $\pi_1(X_n)\cong Q/N(r^{2n+1})$, where $N(S)$ denotes the
normal closure of $S$, and this group is $\Z/(2n+1)\rtimes_{-1}\Z$ by
\cite[Lemma~4.6]{OzSmooth}.  For $n\ge1$ it is not isomorphic to
$\pi_1(X_0)\cong\Z$, since it has torsion, and the groups for different $n$ have
torsion subgroups of different orders.  So the $X_n$ are pairwise
non-homeomorphic.  Finally, suppose that $\Psi^n$, for some $n\ge1$, extends to
a homeomorphism $F$ of $DT^*K$.  Then $F\circ\Psi^{-n}\circ j=j$ on the boundary,
so $F$ together with the identity on the attached handle is a homeomorphism from
$X_n$ to $X_0$, a contradiction.  Hence no positive power of $\Psi$ extends over
$DT^*K$, even as a homeomorphism.
\end{proof}

\begin{remark}\label[remark]{rem:GP}
We expect the Stein structure on $X_{B'}$ given by the Lefschetz fibration
\cite{AO,LP} to be deformation equivalent to $(DT^*K,\omega_{\mathrm{can}})$.
Here and below, deformation equivalence of Stein or Weinstein domains is in the
sense of \cite{GP}, i.e., a Weinstein homotopy of the underlying Weinstein
structures, allowing birth--death of critical points.  Johns' theorem gives a
conformally exact symplectomorphism, which is what we use, but not a Weinstein
deformation equivalence, and we do not claim one here.  Only the smooth identification
of $X_{B'}$ with $DT^*K$ and the identification of the boundary contact
structures enter the proof of \cref{thm:main}.
\end{remark}

\begin{proof}[Proof of \cref{cor:finiteness}]
By \cite[Corollary~1.2]{LiOzbagci}, $(ST^*K,\xi_{\mathrm{can}})$ has a unique minimal weak
symplectic filling up to homeomorphism, and it is planar, supported by
$(\Sigma_0^6,B')$.  By \cref{thm:main}(2), a single Legendrian surgery on it gives
$(Y,\xi)$, which has infinitely many pairwise non-homeomorphic Stein fillings.
Stein fillings are exact, and hence minimal, so $(Y,\xi)$ has infinitely many
minimal symplectic fillings up to homeomorphism.
\end{proof}

The attaching knots are not cotangent-circle fibers: by \cref{lem:pi1}, the curve
$\Psi^{-n}(d')$ represents $r^{-(2n+1)}$ in $\pi_1(DT^*K)\cong Q$, an element of
infinite order, whereas a fiber bounds a disk in $DT^*K$.

The proof of \cref{thm:main}(2) uses the curve $d'$ only as a page curve.  For a
page curve $c$, let $(Y^{B'}_c,\xi^{B'}_c)$ be the contact $3$--manifold supported by the
open book $(\Sigma_0^6,t_cB')$, and let $X^{B'}_n(c)$ be the Lefschetz fibration of
the positive factorization $t_c\,t_{\alpha'}^nB't_{\alpha'}^{-n}$ of $t_cB'$,
$n\ge0$.

\begin{proposition}\label[proposition]{prop:general}
For every page curve $c$, the contact $3$--manifold $(Y^{B'}_c,\xi^{B'}_c)$ is obtained
from $(ST^*K,\xi_{\mathrm{can}})$ by Legendrian surgery along a page realization
of $c$, and $X^{B'}_n(c)\cong DT^*K\cup h_{\Psi^{-n}(c)}$ for every $n\ge0$, where the
$2$--handle is attached with framing one less than the page framing.
\end{proposition}

\begin{proof}
The proof of \cref{thm:main}(2), with $d'$ replaced by $c$, applies verbatim.
\end{proof}

The non-extension statement of \cref{thm:main}(2) has the following
consequence for the mapping class groups of $ST^*K$.

\begin{corollary}\label[corollary]{cor:mcg}
The contactomorphism $\Psi$ has infinite order in the contact mapping class
group of $(ST^*K,\xi_{\mathrm{can}})$, and even in the smooth mapping class group
of $ST^*K$.  More precisely, for every $k\neq0$, the diffeomorphism $\Psi^k$ is
not isotopic, through homeomorphisms, to the restriction of a homeomorphism of
$DT^*K$.
\end{corollary}

\begin{proof}
Suppose that $\Psi^k$, $k\neq0$, is isotopic to the restriction of a
homeomorphism $F$ of $DT^*K$.  Replacing $k$ by $-k$ and $F$ by $F^{-1}$ if
necessary, we may assume that $k\ge1$.  Using the isotopy in a collar of
$ST^*K=\partial DT^*K$, we modify $F$ to a homeomorphism of $DT^*K$ whose
restriction to the boundary is $\Psi^k$, which contradicts \cref{thm:main}(2).
Taking $F$ to be the identity, $\Psi^k$ is not isotopic to the identity for
$k\neq0$, so $\Psi$ has infinite order in the smooth mapping class group of
$ST^*K$.  Since contactomorphisms which are contact isotopic are smoothly
isotopic, $\Psi$ also has infinite order in the contact mapping class group.
\end{proof}

\section{Relation with Baykur's example}\label{sec:baykur}

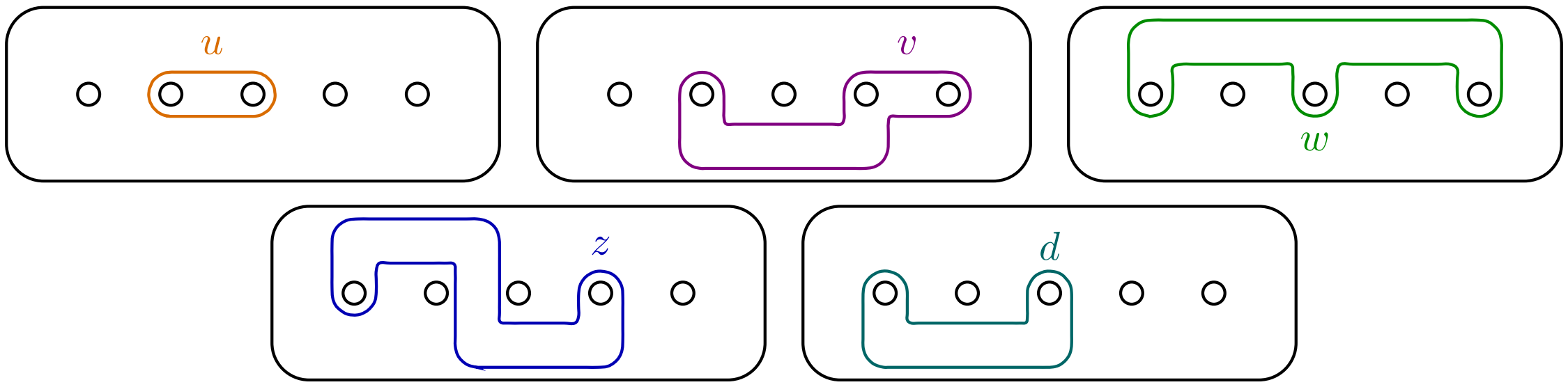
\begin{figure}[ht]
\centering
\makebox[\textwidth][c]{%
\begin{tikzpicture}[scale=0.3]
\begin{scope}[shift={(0.0,0.0)}]
\draw[thick,rounded corners=10pt] (0,-2.75) rectangle (15.600000000000001,2.75);
\draw[thick] (2.60,0.00) circle (0.35);
\draw[thick] (5.20,0.00) circle (0.35);
\draw[thick] (7.80,0.00) circle (0.35);
\draw[thick] (10.40,0.00) circle (0.35);
\draw[thick] (13.00,0.00) circle (0.35);
\draw[orange!85!black,thick] plot[smooth cycle,tension=0.7] coordinates{(5.14,-0.69) (5.05,-0.68) (4.97,-0.66) (4.90,-0.63) (4.83,-0.59) (4.77,-0.55) (4.72,-0.50) (4.67,-0.45) (4.62,-0.39) (4.59,-0.33) (4.55,-0.26) (4.53,-0.19) (4.51,-0.10) (4.50,-0.01) (4.51,0.08) (4.52,0.16) (4.55,0.24) (4.58,0.31) (4.61,0.38) (4.65,0.44) (4.70,0.49) (4.76,0.54) (4.81,0.58) (4.88,0.62) (4.95,0.65) (5.03,0.68) (5.12,0.69) (5.21,0.70) (5.31,0.70) (5.42,0.70) (5.52,0.70) (5.63,0.70) (5.73,0.70) (5.83,0.70) (5.94,0.70) (6.04,0.70) (6.14,0.70) (6.25,0.70) (6.35,0.70) (6.46,0.70) (6.56,0.70) (6.66,0.70) (6.77,0.70) (6.87,0.70) (6.97,0.70) (7.08,0.70) (7.18,0.70) (7.29,0.70) (7.39,0.70) (7.49,0.70) (7.60,0.70) (7.70,0.70) (7.80,0.70) (7.90,0.69) (7.98,0.67) (8.06,0.65) (8.13,0.61) (8.19,0.58) (8.25,0.53) (8.30,0.48) (8.35,0.43) (8.39,0.37) (8.43,0.30) (8.46,0.23) (8.48,0.15) (8.49,0.06) (8.50,-0.03) (8.49,-0.13) (8.47,-0.21) (8.44,-0.28) (8.40,-0.35) (8.36,-0.41) (8.32,-0.47) (8.27,-0.52) (8.21,-0.56) (8.15,-0.60) (8.08,-0.64) (8.01,-0.67) (7.92,-0.69) (7.83,-0.70) (7.73,-0.70) (7.63,-0.70) (7.52,-0.70) (7.42,-0.70) (7.32,-0.70) (7.21,-0.70) (7.11,-0.70) (7.01,-0.70) (6.90,-0.70) (6.80,-0.70) (6.69,-0.70) (6.59,-0.70) (6.49,-0.70) (6.38,-0.70) (6.28,-0.70) (6.18,-0.70) (6.07,-0.70) (5.97,-0.70) (5.86,-0.70) (5.76,-0.70) (5.66,-0.70) (5.55,-0.70) (5.45,-0.70) (5.34,-0.70) (5.24,-0.70) (5.15,-0.70)};
\node[orange!85!black] at (6.50,1.55) {$u$};
\end{scope}
\begin{scope}[shift={(16.8,0.0)}]
\draw[thick,rounded corners=10pt] (0,-2.75) rectangle (15.600000000000001,2.75);
\draw[thick] (2.60,0.00) circle (0.35);
\draw[thick] (5.20,0.00) circle (0.35);
\draw[thick] (7.80,0.00) circle (0.35);
\draw[thick] (10.40,0.00) circle (0.35);
\draw[thick] (13.00,0.00) circle (0.35);
\draw[violet,thick] plot[smooth cycle,tension=0.7] coordinates{(5.17,-2.35) (4.94,-2.30) (4.77,-2.20) (4.63,-2.06) (4.54,-1.88) (4.50,-1.65) (4.50,-1.38) (4.50,-1.11) (4.50,-0.84) (4.50,-0.57) (4.50,-0.30) (4.50,-0.03) (4.53,0.21) (4.62,0.39) (4.76,0.54) (4.93,0.64) (5.15,0.69) (5.40,0.67) (5.58,0.58) (5.73,0.45) (5.84,0.28) (5.89,0.07) (5.90,-0.20) (5.90,-0.47) (5.90,-0.74) (5.97,-0.94) (6.24,-0.95) (6.52,-0.95) (6.80,-0.95) (7.08,-0.95) (7.36,-0.95) (7.64,-0.95) (7.92,-0.95) (8.20,-0.95) (8.48,-0.95) (8.76,-0.95) (9.04,-0.95) (9.32,-0.95) (9.60,-0.95) (9.70,-0.78) (9.70,-0.51) (9.70,-0.24) (9.70,0.03) (9.75,0.26) (9.85,0.43) (9.99,0.57) (10.17,0.66) (10.41,0.70) (10.69,0.70) (10.97,0.70) (11.25,0.70) (11.53,0.70) (11.81,0.70) (12.09,0.70) (12.37,0.70) (12.65,0.70) (12.93,0.70) (13.18,0.67) (13.38,0.59) (13.52,0.46) (13.63,0.29) (13.69,0.08) (13.68,-0.17) (13.60,-0.36) (13.47,-0.51) (13.31,-0.62) (13.10,-0.69) (12.83,-0.70) (12.55,-0.70) (12.27,-0.70) (11.99,-0.70) (11.71,-0.70) (11.43,-0.70) (11.16,-0.71) (11.10,-0.93) (11.10,-1.20) (11.10,-1.47) (11.09,-1.73) (11.03,-1.94) (10.93,-2.11) (10.78,-2.24) (10.59,-2.32) (10.34,-2.35) (10.06,-2.35) (9.78,-2.35) (9.50,-2.35) (9.22,-2.35) (8.94,-2.35) (8.66,-2.35) (8.38,-2.35) (8.10,-2.35) (7.82,-2.35) (7.54,-2.35) (7.25,-2.35) (6.97,-2.35) (6.69,-2.35) (6.41,-2.35) (6.13,-2.35) (5.85,-2.35) (5.57,-2.35) (5.29,-2.35)};
\node[violet] at (11.70,1.55) {$v$};
\end{scope}
\begin{scope}[shift={(33.6,0.0)}]
\draw[thick,rounded corners=10pt] (0,-2.75) rectangle (15.600000000000001,2.75);
\draw[thick] (2.60,0.00) circle (0.35);
\draw[thick] (5.20,0.00) circle (0.35);
\draw[thick] (7.80,0.00) circle (0.35);
\draw[thick] (10.40,0.00) circle (0.35);
\draw[thick] (13.00,0.00) circle (0.35);
\draw[green!55!black,thick] plot[smooth cycle,tension=0.7] coordinates{(2.53,-0.69) (2.25,-0.60) (2.05,-0.43) (1.92,-0.17) (1.90,0.18) (1.90,0.55) (1.90,0.92) (1.90,1.29) (1.90,1.66) (1.97,1.96) (2.14,2.17) (2.38,2.31) (2.70,2.35) (3.00,2.35) (3.32,2.35) (3.62,2.35) (3.94,2.35) (4.24,2.35) (4.55,2.35) (4.87,2.35) (5.17,2.35) (5.48,2.35) (5.79,2.35) (6.10,2.35) (6.41,2.35) (6.72,2.35) (7.03,2.35) (7.34,2.35) (7.65,2.35) (7.96,2.35) (8.27,2.35) (8.58,2.35) (8.89,2.35) (9.20,2.35) (9.51,2.35) (9.82,2.35) (10.13,2.35) (10.44,2.35) (10.75,2.35) (11.06,2.35) (11.37,2.35) (11.68,2.35) (11.99,2.35) (12.30,2.35) (12.61,2.35) (12.92,2.35) (13.24,2.31) (13.47,2.16) (13.63,1.95) (13.70,1.64) (13.70,1.27) (13.70,0.90) (13.70,0.53) (13.70,0.16) (13.67,-0.19) (13.54,-0.44) (13.34,-0.61) (13.04,-0.70) (12.73,-0.64) (12.50,-0.49) (12.35,-0.26) (12.30,0.06) (12.30,0.43) (12.30,0.80) (12.07,0.95) (11.69,0.95) (11.30,0.95) (10.92,0.95) (10.54,0.95) (10.15,0.95) (9.77,0.95) (9.38,0.95) (9.00,0.95) (8.62,0.95) (8.50,0.70) (8.50,0.33) (8.50,-0.04) (8.41,-0.34) (8.24,-0.54) (7.99,-0.67) (7.66,-0.68) (7.39,-0.57) (7.21,-0.37) (7.11,-0.10) (7.10,0.27) (7.10,0.64) (7.04,0.94) (6.66,0.95) (6.28,0.95) (5.89,0.95) (5.51,0.95) (5.13,0.95) (4.74,0.95) (4.36,0.95) (3.97,0.95) (3.59,0.95) (3.30,0.85) (3.30,0.49) (3.30,0.12) (3.26,-0.22) (3.13,-0.46) (2.92,-0.62) (2.61,-0.70)};
\node[green!55!black] at (7.80,-1.50) {$w$};
\end{scope}
\begin{scope}[shift={(8.4,-6.3)}]
\draw[thick,rounded corners=10pt] (0,-2.75) rectangle (15.600000000000001,2.75);
\draw[thick] (2.60,0.00) circle (0.35);
\draw[thick] (5.20,0.00) circle (0.35);
\draw[thick] (7.80,0.00) circle (0.35);
\draw[thick] (10.40,0.00) circle (0.35);
\draw[thick] (13.00,0.00) circle (0.35);
\draw[blue!70!black,thick] plot[smooth cycle,tension=0.7] coordinates{(6.48,-2.35) (6.21,-2.28) (6.01,-2.14) (5.87,-1.95) (5.80,-1.68) (5.80,-1.36) (5.80,-1.03) (5.80,-0.70) (5.80,-0.37) (5.80,-0.04) (5.80,0.30) (5.80,0.63) (5.76,0.92) (5.45,0.95) (5.11,0.95) (4.76,0.95) (4.42,0.95) (4.08,0.95) (3.74,0.95) (3.40,0.95) (3.30,0.71) (3.30,0.38) (3.30,0.05) (3.26,-0.24) (3.13,-0.45) (2.95,-0.60) (2.70,-0.69) (2.40,-0.67) (2.18,-0.55) (2.01,-0.38) (1.92,-0.14) (1.90,0.18) (1.90,0.51) (1.90,0.84) (1.90,1.17) (1.90,1.50) (1.92,1.81) (2.02,2.04) (2.19,2.21) (2.42,2.32) (2.73,2.35) (3.07,2.35) (3.41,2.35) (3.76,2.35) (4.10,2.35) (4.44,2.35) (4.78,2.35) (5.13,2.35) (5.47,2.35) (5.81,2.35) (6.15,2.35) (6.49,2.35) (6.77,2.29) (6.98,2.16) (7.12,1.97) (7.19,1.70) (7.20,1.38) (7.20,1.05) (7.20,0.72) (7.20,0.39) (7.20,0.06) (7.20,-0.28) (7.20,-0.61) (7.22,-0.91) (7.52,-0.95) (7.87,-0.95) (8.21,-0.95) (8.55,-0.95) (8.90,-0.95) (9.24,-0.95) (9.58,-0.95) (9.70,-0.74) (9.70,-0.41) (9.70,-0.08) (9.74,0.22) (9.86,0.44) (10.03,0.59) (10.28,0.69) (10.58,0.67) (10.81,0.57) (10.98,0.39) (11.08,0.16) (11.10,-0.15) (11.10,-0.48) (11.10,-0.81) (11.10,-1.14) (11.10,-1.47) (11.08,-1.78) (10.99,-2.02) (10.83,-2.20) (10.61,-2.32) (10.31,-2.35) (9.96,-2.35) (9.62,-2.35) (9.28,-2.35) (8.94,-2.35) (8.59,-2.35) (8.25,-2.35) (7.91,-2.35) (7.57,-2.35) (7.22,-2.35) (6.88,-2.35) (6.54,-2.35)};
\node[blue!70!black] at (10.40,1.50) {$z$};
\end{scope}
\begin{scope}[shift={(25.2,-6.3)}]
\draw[thick,rounded corners=10pt] (0,-2.75) rectangle (15.600000000000001,2.75);
\draw[thick] (2.60,0.00) circle (0.35);
\draw[thick] (5.20,0.00) circle (0.35);
\draw[thick] (7.80,0.00) circle (0.35);
\draw[thick] (10.40,0.00) circle (0.35);
\draw[thick] (13.00,0.00) circle (0.35);
\draw[teal!80!black,thick] plot[smooth cycle,tension=0.7] coordinates{(2.57,-2.35) (2.39,-2.31) (2.23,-2.24) (2.11,-2.14) (2.01,-2.02) (1.94,-1.87) (1.90,-1.68) (1.90,-1.46) (1.90,-1.24) (1.90,-1.02) (1.90,-0.80) (1.90,-0.58) (1.90,-0.36) (1.90,-0.14) (1.91,0.08) (1.95,0.26) (2.03,0.40) (2.14,0.52) (2.27,0.61) (2.43,0.68) (2.63,0.70) (2.82,0.66) (2.98,0.59) (3.10,0.49) (3.20,0.36) (3.27,0.21) (3.30,0.02) (3.30,-0.20) (3.30,-0.42) (3.30,-0.64) (3.30,-0.85) (3.43,-0.95) (3.66,-0.95) (3.89,-0.95) (4.12,-0.95) (4.35,-0.95) (4.58,-0.95) (4.80,-0.95) (5.03,-0.95) (5.26,-0.95) (5.49,-0.95) (5.72,-0.95) (5.95,-0.95) (6.18,-0.95) (6.40,-0.95) (6.63,-0.95) (6.86,-0.95) (7.06,-0.93) (7.10,-0.75) (7.10,-0.53) (7.10,-0.31) (7.10,-0.09) (7.11,0.13) (7.17,0.29) (7.25,0.43) (7.36,0.55) (7.50,0.63) (7.68,0.69) (7.88,0.69) (8.06,0.65) (8.21,0.57) (8.32,0.46) (8.42,0.33) (8.48,0.17) (8.50,-0.04) (8.50,-0.26) (8.50,-0.48) (8.50,-0.70) (8.50,-0.92) (8.50,-1.14) (8.50,-1.36) (8.50,-1.58) (8.48,-1.79) (8.43,-1.95) (8.34,-2.09) (8.23,-2.20) (8.09,-2.29) (7.91,-2.34) (7.70,-2.35) (7.47,-2.35) (7.24,-2.35) (7.02,-2.35) (6.79,-2.35) (6.56,-2.35) (6.33,-2.35) (6.10,-2.35) (5.87,-2.35) (5.65,-2.35) (5.42,-2.35) (5.19,-2.35) (4.96,-2.35) (4.73,-2.35) (4.50,-2.35) (4.28,-2.35) (4.05,-2.35) (3.82,-2.35) (3.59,-2.35) (3.36,-2.35) (3.13,-2.35) (2.91,-2.35) (2.68,-2.35)};
\node[teal!80!black] at (7.80,1.50) {$d$};
\end{scope}
\end{tikzpicture}}
\caption{Baykur's curves $u,v,w,z,d$ on the page $\Sigma_0^6$.}
\label{fig:baykur-curves}
\end{figure}

In this section we prove \cref{thm:main}(3).  In the generators of
\cref{fig:oo-curves}, Baykur's curves are $u=x_3x_2$, $v=x_5x_4x_2$,
$w=x_5x_4x_3x_2x_1x_2^{-1}x_4^{-1}$, $z=x_4x_2x_1x_2^{-1}$ and $d=x_3x_1$, up to
conjugation \cite{BaykurPlanar}, and $\alpha=x_2x_1$ and $\beta=x_4x_3$ are the
round curves of \cref{fig:oo-curves}.  Baykur's generalized lantern relation
\cite[Lemma~3]{BaykurPlanar} is $Bt_\alpha t_\beta=\Delta$, with
$\Delta=t_{\delta_0}t_{\delta_1}^2t_{\delta_2}^2t_{\delta_3}^2t_{\delta_4}^2t_{\delta_5}$,
and it gives $[B,t_\alpha]=1$.  We write $(Y_B,\xi_B)$ for the contact $3$--manifold
supported by the open book $(\Sigma_0^6,B)$.  The two factorizations
$B=t_ut_vt_wt_z$ and $B'=t_{V_0}t_{V_1}t_{V_2}t_{V_3}$ are not related by a diffeomorphism of the page carrying vanishing cycles to
vanishing cycles: $V_1$ and $V_2$ are disjoint, whereas any two of $u,v,w,z$
intersect.  They become so related after three Hurwitz moves.  All that is
needed are the following identities between curves, where $\beta'$ is the curve
of \cref{fig:oo-curves}.

\begin{lemma}\label[lemma]{lem:curve-identities}
Let
\[
g=\sigma_{\alpha'}\,\sigma_3^{-1}\sigma_2^{-1}\sigma_4,\qquad
\sigma_{\alpha'}=\sigma_1\sigma_2\,\sigma_3^{-1}\,\sigma_2^{-1}\sigma_1^{-1},
\]
in functional order, with $\sigma_i$ as in \eqref{eq:artin}.  Here
$\sigma_{\alpha'}$ is the counterclockwise half twist supported in the disk bounded by
$\alpha'=\sigma_1\sigma_2(\beta)$, interchanging $\delta_1$ and $\delta_4$, since it
is conjugate to $\sigma_3^{-1}$.  Then $g$
acts on $\pi_1(\Sigma_0^6,\ast)$ by
\begin{gather*}
x_1\mapsto x_4,\qquad x_2\mapsto x_4x_1x_4^{-1},\qquad
x_3\mapsto x_4x_2x_4^{-1},\qquad x_5\mapsto x_4x_3x_4^{-1},\\
x_4\mapsto x_4x_3^{-1}x_4^{-1}\,x_5\,x_4x_3x_4^{-1}.
\end{gather*}
Moreover, $g$ is the identity on $\delta_0$, permutes the holes by
$\delta_1\mapsto\delta_4$, $\delta_2\mapsto\delta_1$,
$\delta_3\mapsto\delta_2$, $\delta_4\mapsto\delta_5$,
$\delta_5\mapsto\delta_3$, and, up to isotopy,
\begin{gather*}
g(d)=d',\qquad g(\alpha)=\alpha',\qquad g(\beta)=\beta',\\
g(u)=V_1,\qquad g(v)=V_3,\qquad g\bigl(t_v(z)\bigr)=V_2,\qquad g\bigl(t_ut_v(w)\bigr)=V_0 .
\end{gather*}
\end{lemma}

\begin{proof}
As explained in \cref{sec:openbook}, it suffices to compare conjugacy classes in
$\pi_1(\Sigma_0^6,\ast)=F(x_1,\dots,x_5)$, up to inversion.  Since $g$ is
the identity on $\delta_0\ni\ast$, it induces an automorphism of this
group, and since only identities between curves are asserted, nothing is
lost in passing to $\pi_1$: the action of the mapping class group on
$\pi_1(\Sigma_0^6,\ast)$ fails to detect only twists along the inner
boundary components, which play no role here.

We first compute the action of $g$.  Composing the formulas
\eqref{eq:artin}--\eqref{eq:artininv} one factor at a time, from the right, we
obtain the following automorphisms, where generators not listed are fixed:
\begin{align*}
\sigma_1^{-1}\colon&\quad x_1\mapsto x_2,\quad x_2\mapsto x_2x_1x_2^{-1},\\
\sigma_2^{-1}\sigma_1^{-1}\colon&\quad x_1\mapsto x_3,\quad x_2\mapsto x_3x_1x_3^{-1},\quad x_3\mapsto x_3x_2x_3^{-1},\\
\sigma_3^{-1}\sigma_2^{-1}\sigma_1^{-1}\colon&\quad x_1\mapsto x_4,\quad x_2\mapsto x_4x_1x_4^{-1},\quad x_3\mapsto x_4x_2x_4^{-1},\quad x_4\mapsto x_4x_3x_4^{-1},\\
\sigma_2\sigma_3^{-1}\sigma_2^{-1}\sigma_1^{-1}\colon&\quad x_1\mapsto x_4,\quad x_2\mapsto x_4x_1x_4^{-1},\quad x_3\mapsto x_4x_2^{-1}x_3x_2x_4^{-1},\quad x_4\mapsto x_4x_2x_4^{-1},\\
\sigma_{\alpha'}\colon&\quad x_1\mapsto x_4,\quad x_2\mapsto x_4x_1^{-1}x_2x_1x_4^{-1},\quad x_3\mapsto x_4x_1^{-1}x_3x_1x_4^{-1},\quad x_4\mapsto x_4x_1x_4^{-1}.
\end{align*}
For example, in the second line, $\sigma_2^{-1}\sigma_1^{-1}(x_2)$ is obtained by
applying $\sigma_2^{-1}$ to $\sigma_1^{-1}(x_2)=x_2x_1x_2^{-1}$, which gives
$x_3x_1x_3^{-1}$.  Similarly,
\begin{align*}
\sigma_4\colon&\quad x_4\mapsto x_4^{-1}x_5x_4,\quad x_5\mapsto x_4,\\
\sigma_2^{-1}\sigma_4\colon&\quad x_2\mapsto x_3,\quad x_3\mapsto x_3x_2x_3^{-1},\quad x_4\mapsto x_4^{-1}x_5x_4,\quad x_5\mapsto x_4,\\
\sigma_3^{-1}\sigma_2^{-1}\sigma_4\colon&\quad x_2\mapsto x_4,\quad x_3\mapsto x_4x_2x_4^{-1},\quad x_4\mapsto x_4x_3^{-1}x_4^{-1}x_5x_4x_3x_4^{-1},\quad x_5\mapsto x_4x_3x_4^{-1}.
\end{align*}
Finally $g=\sigma_{\alpha'}\circ(\sigma_3^{-1}\sigma_2^{-1}\sigma_4)$, so
$g(x_i)$ is obtained by applying $\sigma_{\alpha'}$ to the image of $x_i$ in the
last line.  This gives $g(x_1)=\sigma_{\alpha'}(x_1)=x_4$,
$g(x_2)=\sigma_{\alpha'}(x_4)=x_4x_1x_4^{-1}$,
\[
g(x_3)=\sigma_{\alpha'}(x_4x_2x_4^{-1})
=(x_4x_1x_4^{-1})(x_4x_1^{-1}x_2x_1x_4^{-1})(x_4x_1^{-1}x_4^{-1})=x_4x_2x_4^{-1},
\]
and in the same way, since $\sigma_{\alpha'}(x_3)=x_4x_1^{-1}x_3x_1x_4^{-1}$,
\[
g(x_5)=\sigma_{\alpha'}(x_4x_3x_4^{-1})
=(x_4x_1x_4^{-1})(x_4x_1^{-1}x_3x_1x_4^{-1})(x_4x_1^{-1}x_4^{-1})=x_4x_3x_4^{-1}.
\]
Writing $x_4x_3^{-1}x_4^{-1}x_5x_4x_3x_4^{-1}=(x_4x_3x_4^{-1})^{-1}\,x_5\,(x_4x_3x_4^{-1})$
and using $\sigma_{\alpha'}(x_5)=x_5$, we get
$g(x_4)=g(x_5)^{-1}\,x_5\,g(x_5)=x_4x_3^{-1}x_4^{-1}x_5x_4x_3x_4^{-1}$, as claimed.

The classes of the curves involved are
\begin{gather*}
\alpha=x_2x_1,\qquad \beta=x_4x_3,\qquad u=x_3x_2,\qquad v=x_5x_4x_2,\\
t_v(z)=x_5^{-1}x_1x_5\,x_4,\qquad t_ut_v(w)=x_1x_5x_3,\qquad
\alpha'=x_1x_4,\qquad \beta'=x_4^{-1}x_5x_4\,x_3x_2x_3^{-1},\\
V_0=x_4x_3x_2,\qquad V_1=x_2x_1,\qquad V_2=x_5x_4,\qquad V_3=x_5x_4x_3x_1x_4^{-1}.
\end{gather*}
The curves $t_v(z)$ and
$t_ut_v(w)$, obtained from the curves of \cref{fig:baykur-curves} by
twisting $z$ along $v$, and $w$ along $v$ and then along $u$, are shown in
\cref{fig:twisted-curves}, from which their classes are read off.  These two
classes can also be computed algebraically, with the description of round twists
in \cref{sec:openbook}.  Thus $t_u$
conjugates $x_2,x_3$ by $u=x_3x_2$.  Let $\rho$ be the boundary of the disk
containing $\delta_3,\delta_4,\delta_5$, represented by $x_5x_4x_3$.  Then
$\sigma_2(\rho)=x_5x_4x_2=v$, so $t_v=\sigma_2\,t_\rho\,\sigma_2^{-1}$.  Composing
$\sigma_2^{-1}$, then $t_\rho$, which conjugates $x_3,x_4,x_5$ by $\rho$, and
then $\sigma_2$, one finds that $t_v$ fixes $x_1$, conjugates $x_2,x_4,x_5$ by
$v=x_5x_4x_2$, and maps $x_3$ to
$(x_5x_4x_2x_4^{-1}x_5^{-1}x_2^{-1})\,x_3\,(x_5x_4x_2x_4^{-1}x_5^{-1}x_2^{-1})^{-1}$.
Substituting into $z=x_4x_2x_1x_2^{-1}$ and $w=x_5x_4x_3x_2x_1x_2^{-1}x_4^{-1}$ and
freely reducing gives
$t_v(w)=x_5x_4x_3x_2x_1x_5x_2^{-1}x_4^{-1}x_5^{-1}$, and then applying $t_u$,
which conjugates $x_2,x_3$ by $u=x_3x_2$ and fixes the other generators, gives
\begin{gather*}
t_v(z)=(x_5x_4x_2)\,x_5^{-1}x_1x_5x_4\,(x_5x_4x_2)^{-1},\\
t_ut_v(w)=(x_5x_4x_3x_2)\,x_1x_5x_3\,(x_5x_4x_3x_2)^{-1}.
\end{gather*}

\begin{figure}[ht]
\centering
\begin{tikzpicture}[scale=0.42]
\begin{scope}[shift={(0,0)}]
\draw[thick,rounded corners=10pt] (-0.3,-2.75) rectangle (12.6,2.75);
\draw[thick] (2.00,0.00) circle (0.35);
\draw[thick] (4.00,0.00) circle (0.35);
\draw[thick] (6.00,0.00) circle (0.35);
\draw[thick] (8.00,0.00) circle (0.35);
\draw[thick] (10.00,0.00) circle (0.35);
\draw[blue!70!black,thick] plot[smooth cycle,tension=0.7] coordinates{(1.97,-2.37) (1.59,-2.24) (1.35,-1.96) (1.28,-1.49) (1.28,-0.94) (1.28,-0.39) (1.30,0.14) (1.48,0.49) (1.80,0.69) (2.23,0.68) (2.54,0.48) (2.71,0.11) (2.72,-0.43) (2.77,-0.92) (3.26,-0.93) (3.75,-0.93) (4.25,-0.93) (4.74,-0.93) (5.24,-0.93) (5.73,-0.93) (6.23,-0.93) (6.72,-0.93) (7.22,-0.93) (7.71,-0.93) (8.21,-0.93) (8.70,-0.93) (9.20,-0.93) (9.69,-0.93) (10.19,-0.93) (10.62,-0.86) (10.63,-0.32) (10.63,0.24) (10.63,0.79) (10.26,0.93) (9.77,0.93) (9.27,0.93) (8.78,0.92) (8.72,0.45) (8.71,-0.10) (8.55,-0.47) (8.24,-0.68) (7.81,-0.69) (7.49,-0.50) (7.30,-0.16) (7.28,0.37) (7.28,0.92) (7.28,1.47) (7.34,1.94) (7.58,2.23) (7.95,2.37) (8.44,2.37) (8.94,2.37) (9.43,2.37) (9.93,2.37) (10.42,2.37) (10.92,2.37) (11.41,2.37) (11.78,2.23) (12.01,1.93) (12.07,1.45) (12.07,0.90) (12.07,0.35) (12.07,-0.21) (12.07,-0.76) (12.07,-1.31) (12.04,-1.83) (11.85,-2.17) (11.52,-2.35) (11.10,-2.37) (10.70,-2.37) (10.30,-2.37) (9.90,-2.37) (9.49,-2.37) (9.09,-2.37) (8.69,-2.37) (8.29,-2.37) (7.88,-2.37) (7.48,-2.37) (7.08,-2.37) (6.68,-2.37) (6.27,-2.37) (5.87,-2.37) (5.47,-2.37) (5.07,-2.37) (4.67,-2.37) (4.26,-2.37) (3.86,-2.37) (3.46,-2.37) (3.06,-2.37) (2.65,-2.37) (2.25,-2.37)};
\node[blue!70!black] at (1.90,1.50) {$t_v(z)$};
\end{scope}
\begin{scope}[shift={(14.0,0)}]
\draw[thick,rounded corners=10pt] (-0.3,-2.75) rectangle (12.6,2.75);
\draw[thick] (2.00,0.00) circle (0.35);
\draw[thick] (4.00,0.00) circle (0.35);
\draw[thick] (6.00,0.00) circle (0.35);
\draw[thick] (8.00,0.00) circle (0.35);
\draw[thick] (10.00,0.00) circle (0.35);
\draw[green!55!black,thick] plot[smooth cycle,tension=0.7] coordinates{(1.97,-2.37) (1.70,-2.30) (1.49,-2.16) (1.35,-1.95) (1.28,-1.66) (1.28,-1.29) (1.28,-0.92) (1.28,-0.55) (1.28,-0.18) (1.30,0.17) (1.41,0.42) (1.59,0.59) (1.83,0.70) (2.13,0.71) (2.37,0.61) (2.56,0.45) (2.68,0.22) (2.72,-0.12) (2.72,-0.49) (2.72,-0.85) (2.99,-0.93) (3.32,-0.93) (3.65,-0.93) (3.98,-0.93) (4.32,-0.93) (4.65,-0.93) (4.98,-0.93) (5.27,-0.87) (5.28,-0.52) (5.28,-0.15) (5.31,0.19) (5.43,0.43) (5.61,0.60) (5.85,0.70) (6.15,0.70) (6.39,0.60) (6.57,0.43) (6.69,0.19) (6.72,-0.15) (6.72,-0.52) (6.73,-0.87) (7.01,-0.93) (7.34,-0.93) (7.68,-0.93) (8.01,-0.93) (8.34,-0.93) (8.68,-0.93) (9.01,-0.93) (9.27,-0.85) (9.28,-0.49) (9.28,-0.12) (9.31,0.21) (9.44,0.45) (9.62,0.61) (9.87,0.71) (10.17,0.70) (10.41,0.59) (10.58,0.42) (10.70,0.17) (10.72,-0.18) (10.72,-0.55) (10.72,-0.92) (10.72,-1.29) (10.72,-1.66) (10.65,-1.95) (10.50,-2.16) (10.30,-2.30) (10.02,-2.37) (9.69,-2.37) (9.36,-2.37) (9.03,-2.37) (8.70,-2.37) (8.36,-2.37) (8.03,-2.37) (7.70,-2.37) (7.36,-2.37) (7.03,-2.37) (6.70,-2.37) (6.36,-2.37) (6.03,-2.37) (5.70,-2.37) (5.36,-2.37) (5.03,-2.37) (4.71,-2.37) (4.38,-2.37) (4.05,-2.37) (3.71,-2.37) (3.38,-2.37) (3.05,-2.37) (2.72,-2.37) (2.38,-2.37) (2.05,-2.37)};
\node[green!55!black] at (6.00,1.55) {$t_ut_v(w)$};
\end{scope}
\end{tikzpicture}
\caption{The curves $t_v(z)$ and $t_ut_v(w)$ of \cref{lem:curve-identities}.}
\label{fig:twisted-curves}
\end{figure}
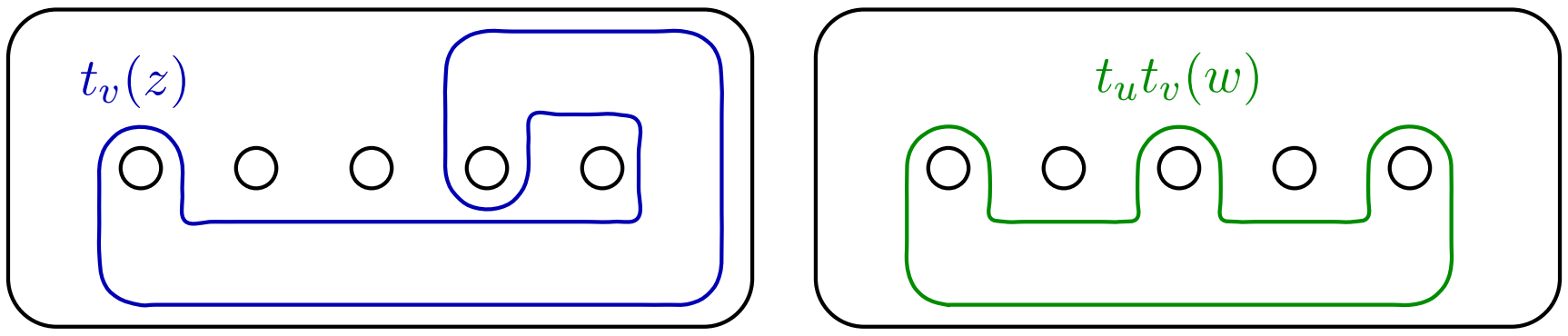
  Applying the
formulas for $g$ and reducing,
\begin{align*}
g(d)&=g(x_3)\,g(x_1)=(x_4x_2x_4^{-1})(x_4)=x_4x_2=d',\\
g(\alpha)&=(x_4x_1x_4^{-1})(x_4)=x_4x_1\sim x_1x_4=\alpha',\\
g(\beta)&=x_4x_3^{-1}x_4^{-1}x_5x_4x_3\,x_2x_4^{-1}
\sim x_3^{-1}x_4^{-1}x_5x_4x_3x_2\sim\beta',\\
g(u)&=x_4x_2x_1x_4^{-1}\sim x_2x_1=V_1,\\
g(v)&=x_5x_4x_3x_1x_4^{-1}=V_3,\\
g\bigl(t_v(z)\bigr)&=x_4x_3^{-1}\,x_5x_4\,x_3x_4^{-1}\sim x_5x_4=V_2,\\
g\bigl(t_ut_v(w)\bigr)&=(x_4)(x_4x_3x_4^{-1})(x_4x_2x_4^{-1})
=x_4\,x_4x_3x_2\,x_4^{-1}\sim V_0 ,
\end{align*}
where $\sim$ denotes conjugacy.  The permutation of the holes is read off
from the images of $x_1,\dots,x_5$, each a conjugate of a generator, and $g$
is the identity near $\delta_0$ since its support can be chosen away from it.
\end{proof}

\begin{proposition}\label[proposition]{prop:baykur-oba-conjugacy}
Three Hurwitz moves transform Baykur's factorization into
\begin{equation}\label{eq:baykur-hurwitz}
B=t_ut_vt_wt_z=t_{t_ut_v(w)}\;t_u\;t_{t_v(z)}\;t_v ,
\end{equation}
and the diffeomorphism $g$ of \cref{lem:curve-identities} carries the four
curves on the right to the curves $V_0,V_1,V_2,V_3$ of \cref{fig:oo-curves}.  In
particular $gBg^{-1}=B'$ and $(Y_B,\xi_B)\cong(ST^*K,\xi_{\mathrm{can}})$.
\end{proposition}

\begin{proof}
The identity \eqref{eq:baykur-hurwitz} is obtained from $t_ut_vt_wt_z$ by the
elementary moves $t_at_b=t_{t_a(b)}t_a$, applied to the second, third and
first pairs of factors in turn:
\[
(u,v,w,z)\to\bigl(u,t_v(w),v,z\bigr)\to\bigl(u,t_v(w),t_v(z),v\bigr)
\to\bigl(t_ut_v(w),u,t_v(z),v\bigr).
\]
Since $g\,t_c\,g^{-1}=t_{g(c)}$ for every simple closed curve $c$,
\cref{lem:curve-identities} gives
\[
gBg^{-1}=t_{g(t_ut_v(w))}\,t_{g(u)}\,t_{g(t_v(z))}\,t_{g(v)}
=t_{V_0}t_{V_1}t_{V_2}t_{V_3}=B' .
\]
Conjugate monodromies give isomorphic open books \cite[\S2]{EtnyreLectures}:
the conjugating diffeomorphism need only preserve the orientation of the page
and map $\partial\Sigma_0^6$ to itself.  Here, for a representative of $B$, denoted by the
same letter, the map $(x,t)\mapsto(g(x),t)$ from the mapping torus of $B$ to
that of $gBg^{-1}$, which represents $B'$, extended by $g|_{\partial\Sigma_0^6}\times\mathrm{id}$
on $\partial\Sigma_0^6\times D^2$, is an isomorphism of open books
$(\Sigma_0^6,B)\to(\Sigma_0^6,B')$ which permutes the binding
components.  Hence $(Y_B,\xi_B)$ is
contactomorphic to the contact manifold supported by $(\Sigma_0^6,B')$,
which is $(ST^*K,\xi_{\mathrm{can}})$ by \cite[Theorem~3.7]{ObaOzbagci}.
\end{proof}

\begin{remark}\label[remark]{rem:transport}
\Cref{prop:lantern} is Baykur's generalized lantern relation
$Bt_\alpha t_\beta=\Delta$, with
$\Delta=t_{\delta_0}t_{\delta_1}^2t_{\delta_2}^2t_{\delta_3}^2t_{\delta_4}^2t_{\delta_5}$
\cite[Lemma~3]{BaykurPlanar}, transported by $g$.  Indeed, by
\cref{lem:curve-identities}, $g(\alpha)=\alpha'$ and $g(\beta)=\beta'$, and $g$ carries
$\delta_5$, the hole with exponent one in $\Delta$, to $\delta_3$, so
$g\Delta g^{-1}=\Delta'$.  Together with \cref{prop:baykur-oba-conjugacy} this gives
$B't_{\alpha'}t_{\beta'}=g(Bt_\alpha t_\beta)g^{-1}=g\Delta g^{-1}=\Delta'$.
\end{remark}

\begin{remark}\label[remark]{rem:search}
The diffeomorphism $g$ and the Hurwitz moves were found by a computer search
in the action on $\pi_1(\Sigma_0^6,\ast)$ of the group $\widehat\Gamma$ of
\cite[\S4.1]{OzSmooth}, which contains $g$.
The verification above is by hand.  In hindsight, $g$ is explained by
\cref{prop:lantern}, which is proved independently: any diffeomorphism
carrying $(\alpha,\beta)$ to $(\alpha',\beta')$ and $\delta_5$ to $\delta_3$ conjugates
$B=\Delta t_\beta^{-1}t_\alpha^{-1}$ to $\Delta't_{\beta'}^{-1}t_{\alpha'}^{-1}=B'$.
\end{remark}

\begin{proof}[Proof of \cref{thm:main}(3)]
By \cref{prop:baykur-oba-conjugacy} and \cref{lem:curve-identities}, the
diffeomorphism $g$ satisfies $gBg^{-1}=B'$, $g(d)=d'$ and $g(\alpha)=\alpha'$.
Since $gt_cg^{-1}=t_{g(c)}$ for every simple closed curve $c$, conjugation by $g$
carries Baykur's monodromy $t_dB$ to $t_{d'}B'$.  As in the proof of
\cref{prop:baykur-oba-conjugacy}, conjugate monodromies give isomorphic open
books, so Baykur's contact $3$--manifold is contactomorphic to the one supported
by $(\Sigma_0^6,t_{d'}B')$, which is $(Y,\xi)$.

For the fillings we compare ordered vanishing cycles.  A Lefschetz fibration over
the disk is determined up to diffeomorphism by its fiber and the isotopy classes
of its ordered vanishing cycles, and Hurwitz moves do not change it up to
diffeomorphism \cite[\S8.2]{GS}.  Baykur's filling $X^{\mathrm B}_n$ has the ordered
vanishing cycles $d,\,t_\alpha^n(u),\,t_\alpha^n(v),\,t_\alpha^n(w),\,t_\alpha^n(z)$.
Keeping $d$ fixed, we apply to the last four the three Hurwitz moves of
\cref{prop:baykur-oba-conjugacy}, transported by $t_\alpha^n$, which gives
$d,\,t_\alpha^n(t_ut_v(w)),\,t_\alpha^n(u),\,t_\alpha^n(t_v(z)),\,t_\alpha^n(v)$.
Applying $g$ and using $gt_\alpha g^{-1}=t_{\alpha'}$ and \cref{lem:curve-identities},
this list becomes $d',\,t_{\alpha'}^n(V_0),\,t_{\alpha'}^n(V_1),\,t_{\alpha'}^n(V_2),\,t_{\alpha'}^n(V_3)$,
which are the ordered vanishing cycles of the factorization
$P'_n=t_{d'}\,t_{\alpha'}^nB't_{\alpha'}^{-n}$.  Hence $X^{\mathrm B}_n$ is
diffeomorphic to the Lefschetz fibration of $P'_n$, which is diffeomorphic to
$X_n$ by the proof of \cref{thm:main}(2).
\end{proof}

\begin{proof}[Proof of \cref{cor:qhb}]
As in the proof of \cref{thm:main}(3), conjugation by $g$ carries the open book
$(\Sigma_0^6,t_{d_m}B)$ to $(\Sigma_0^6,t_{g(d_m)}B')$.  Moreover, the three
Hurwitz moves of \cref{prop:baykur-oba-conjugacy}, transported by $t_\alpha^n$ and
applied with $t_{d_m}$ held fixed, followed by $g$, carry the ordered vanishing
cycles of the factorization $t_{d_m}t_\alpha^nBt_\alpha^{-n}$, whose Lefschetz
fibration is $X_{m,n}$, to those of $t_{g(d_m)}t_{\alpha'}^nB't_{\alpha'}^{-n}$.  By \cref{prop:general} with
$c=g(d_m)$, the contact $3$--manifold $(Y_m,\xi_m)$ is obtained from
$(ST^*K,\xi_{\mathrm{can}})$ by Legendrian surgery along a page realization of
$g(d_m)$, and $X_{m,n}\cong DT^*K\cup h_{\Psi^{-n}(g(d_m))}$.  The same argument
applies to $d'_m$.  By \cite[Theorem~1.1]{OzSmooth} these fillings are rational
homology balls, pairwise non-homeomorphic for fixed $m$, and
$(ST^*K,\xi_{\mathrm{can}})$ has a unique minimal weak symplectic filling up to
homeomorphism by \cite[Corollary~1.2]{LiOzbagci}.
\end{proof}

\Cref{cor:qhb} also explains why these fillings are rational homology balls.
The disk cotangent bundle $DT^*K$ has the rational homology of a circle, since
$H_1(DT^*K;\Z)\cong H_1(K;\Z)\cong\Z\oplus\Z/2$ and $H_2(DT^*K;\Z)=0$.  Under the
diffeomorphism induced by $g$, the attaching circle of $X_{m,n}$ corresponds to
$t_\alpha^{-n}(d_m)$, whose image in $\pi_1(K)=\langle r,b\mid
brb^{-1}=r^{-1}\rangle$, under the identification of \cite[\S4.2]{OzSmooth}, is
$r^{-(4mn+1)}b^2$ \cite[Lemma~5.2]{OzSmooth}.  Its image
in $H_1(DT^*K;\mathbb Q)\cong\mathbb Q$, generated by the image of $b$, is
$2\neq0$, and for $d'_m$ it is $-2$ \cite[Lemma~5.2]{OzSmooth}.  So the
$2$--handle kills the rational first homology and creates no rational second
homology.

\section{A Stein diagram}\label{sec:diagrams}

\begin{figure}[htbp]
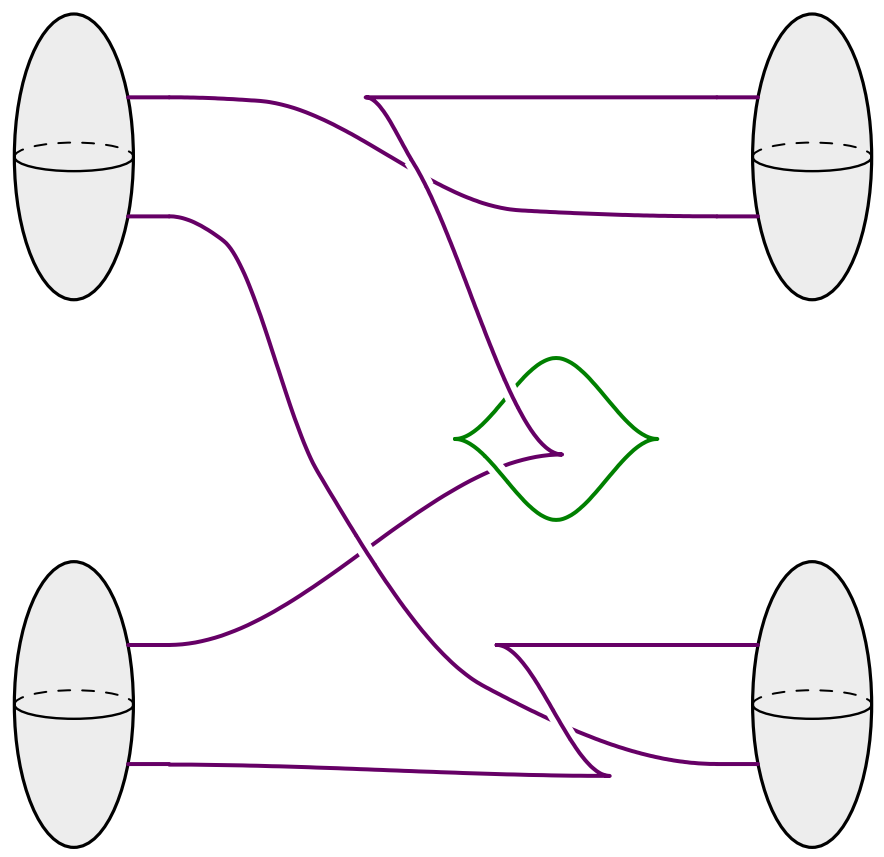

\centering
\scalebox{0.75}{
}
\caption{A Stein structure on $W_0$.  Removing the green knot gives Gompf's
Stein diagram of $DT^*K$.}
\label{fig:W0-stein}
\end{figure}

The fillings $X_n$ of \cref{thm:main} are obtained from $DT^*K$ by attaching a
single Weinstein $2$--handle.  \Cref{fig:W0-stein} shows a Stein domain of the same
shape, which is diffeomorphic to the filling $W_0$ of \cite[Example~4.7]{OzSmooth}
of $(Y,\xi)$.  For $W_0$ we use Baykur's description, which is equivalent to ours by
\cref{thm:main}(3).

The purple Legendrian knot, with $\mathrm{tb}=-3$, together with the two pairs of
balls, is Gompf's Stein diagram \cite{Gompf,OzArch} of $DT^*K$, which is
deformation equivalent to $(DT^*K,\omega_{\mathrm{can}})$ by \cite[Theorem~7.3]{Acu}.
The green curve is a Legendrian unknot with $\mathrm{tb}=-1$ and $\mathrm{rot}=0$, a
meridian of the purple knot.  After the purple handle is attached, it becomes a
Legendrian knot in $ST^*K$ representing the belt circle of the $2$--handle of
$DT^*K$, which is smoothly a circle fiber of $ST^*K\to K$, and attaching a Weinstein
$2$--handle along it gives a Stein domain diffeomorphic to $W_0$ by
\cite[Proposition~4.8]{OzSmooth}.  The generator of $H_2(W_0;\Z)\cong\Z$ is represented by
the cellular $2$--cycle corresponding to the $2$--handle along the green knot.  Its
attaching circle is nullhomologous in the $1$--handlebody, since it runs over no
$1$--handle, and its rotation number is zero, so $c_1$ evaluates to zero on it, and
$c_1=0$ in $H^2(W_0;\mathbb Q)$.  The same holds for the Stein
structure of the Lefschetz fibration $W_0$, whose generator is represented by the
difference of the $2$--handles along two parallel copies of $d$.  Since $\chi(W_0)=1$
and $\sigma(W_0)=-1$, both induced contact structures on $Y$ have $d_3$--invariant
$\frac14(c_1^2-2\chi-3\sigma)=\frac14$.

\begin{conjecture}\label[conjecture]{conj:W0-stein}
The Stein domain of \cref{fig:W0-stein} is deformation equivalent to the Stein
filling $W_0$ of \cite[Example~4.7]{OzSmooth}.  In particular, $(Y,\xi)$ is obtained
from $(ST^*K,\xi_{\mathrm{can}})$ by Legendrian surgery along a Legendrian
realization of a circle fiber of $ST^*K\to K$.
\end{conjecture}

Similarly, a Legendrian realization of the Kirby diagram of $X_0$ on the left of
\cite[Figure~4]{OzSmooth}, with one $1$--handle and a $2$--handle along a Legendrian
trefoil with $\mathrm{tb}=-7$ and $\mathrm{rot}=0$, gives a Stein structure on $X_0$
with $c_1=0$ and $d_3=\frac14$, as for the Lefschetz fibration of $t_{d'}B'$, and we
expect the two Stein structures to be deformation equivalent.

\section{Rational homology ball fillings from Luttinger twists}\label{sec:luttinger}

Luttinger surgery along a Lagrangian Klein bottle, introduced by Nemirovski
\cite{N02,N09b,N09}, is a four-dimensional cut-and-paste operation: a Weinstein neighborhood of
the Klein bottle, which is a copy of $DT^*K$, is cut out and reglued along its boundary $ST^*K$ by
a diffeomorphism of $ST^*K$, called a Luttinger twist.  We give
explicit strict contactomorphisms of $(ST^*K,\xi_{\mathrm{can}})$ isotopic to the Luttinger
twists, and we use them to recover, without reference to an open book or a Lefschetz
fibration, the phenomenon of \cref{cor:qhb}: for every Legendrian knot in
$(ST^*K,\xi_{\mathrm{can}})$ in a suitable homotopy class, the resulting contact $3$--manifold has
infinitely many pairwise non-homeomorphic rational homology ball Stein fillings
(\cref{thm:qhb-canonical}).

\subsection{The Luttinger twists}\label{subsec:twists}

We use the flat model of \cite[\S2.1]{N09b}.  Let $T=\mathbb R^2/2\pi\Z^2$ with
coordinates $(\varphi,\psi)$, and let $(r,\theta)$ be the dual
polar coordinates on the fibers of $T^*T$, so that a covector is
$r\cos\theta\,d\varphi+r\sin\theta\,d\psi$.  Then
\[
\lambda_{\mathrm{can}}=r(\cos\theta\,d\varphi+\sin\theta\,d\psi).
\]
The Klein bottle is $K=T/\bigl((\varphi,\psi)\sim(\varphi+\pi,-\psi)\bigr)$.  The derivative of
$(\varphi,\psi)\mapsto(\varphi+\pi,-\psi)$ maps the covector with polar coordinates $(r,\theta)$ to the one with
$(r,-\theta)$, so $T^*K$ is the quotient of
$T^*T$ by $(\varphi,\psi,r,\theta)\sim(\varphi+\pi,-\psi,r,-\theta)$, and $\lambda_{\mathrm{can}}$ descends to it.  For the flat
metric, $ST^*K=\{r=1\}$, and
\[
\alpha=\lambda_{\mathrm{can}}|_{ST^*K}=\cos\theta\,d\varphi+\sin\theta\,d\psi
\]
is a contact form with kernel $\xi_{\mathrm{can}}$.  Note that
$d\alpha=-\sin\theta\,d\theta\wedge d\varphi+\cos\theta\,d\theta\wedge d\psi$, so every term of $d\alpha$ contains $d\theta$.

The projection to $\varphi\in\mathbb R/\pi\Z$ exhibits $ST^*K$ as a torus bundle over the circle, with
fiber the torus with coordinates $(\psi,\theta)$ and monodromy $-I$.  Let $s$ and $t$ be the loops
through $(\varphi,\psi,\theta)=(0,0,0)$ in the $\psi$-- and $\theta$--directions, and let $c$ be the path
$\varphi\mapsto(\varphi,0,0)$, $0\le \varphi\le\pi$, which is a closed loop since $(\pi,0,0)\sim(0,0,0)$, and which maps
once around the base circle of the torus bundle.  Then
\[
H=\pi_1(ST^*K)=\Gamma\rtimes_{-1}\Z=\langle s,t,c\mid [s,t]=1,\ csc^{-1}=s^{-1},\ ctc^{-1}=t^{-1}\rangle,
\qquad \Gamma=\langle s,t\rangle\cong\Z^2 .
\]
We write the elements of $\Gamma$, the fundamental group of a torus fiber, additively.  Since conjugation by $c$ acts on $\Gamma$ by $-1$, a free homotopy class
of loops in $\Gamma$ is an element of $\Gamma$ up to sign.  The loop $t$ is a fiber of $ST^*K\to K$, and the
map $H\to\pi_1(DT^*K)=\pi_1(K)$ kills $t$ and maps $s$ and $c$ to the classes $r$ and $b$ of the loop
$\{\varphi=0\}$ and of the image of $c$, so that $\pi_1(K)=\langle r,b\mid brb^{-1}=r^{-1}\rangle$, with the same presentation as in \cref{lem:pi1}, and its
kernel is $\langle t\rangle$.  Since $c$ acts on $\Gamma$ by $-1$, the element $c^2$ is central in $H$, and so is its
image $b^2$ in $\pi_1(K)$.

For $k\in\Z$, the Luttinger twist of \cite[\S2.1]{N09b} is the diffeomorphism
\[
f_{0,k}(\varphi,\psi,\theta)=(\varphi,\psi+k\theta,\theta)
\]
of $ST^*K$, which is compatible with the relation above.  It fixes the base point $(0,0,0)$
and the loops $s$ and $c$, and it maps $t$ to the loop $\theta\mapsto(0,k\theta,\theta)$, which is homotopic to
$t+ks$.  So it acts on $H$ by $s\mapsto s$, $t\mapsto t+ks$, $c\mapsto c$.

\subsection{Luttinger twists as contactomorphisms}\label{subsec:strict}

The construction in \cite[\S6]{N09}, where the gluing maps are built from two maps that are
contactomorphisms for contact forms whose kernels are isotopic to $\xi_{\mathrm{can}}$
(\cref{rem:nemirovski-maps}), also implies that the Luttinger twists are, up to isotopy,
contactomorphisms.  We give explicit strict
representatives for the contact form $\alpha$.  For $k\neq0$ the twists $f_{0,k}$ are not contactomorphisms of $(ST^*K,\xi_{\mathrm{can}})$, since
$f_{0,k}^*\alpha=\cos\theta\,d\varphi+\sin\theta\,(d\psi+k\,d\theta)=\alpha+k\sin\theta\,d\theta$, so that $\alpha(\partial_\theta)=0$ while
$(f_{0,k}^*\alpha)(\partial_\theta)=k\sin\theta$.  However, the difference is exact,
\[
f_{0,k}^*\alpha=\alpha+dh_k,\qquad h_k=-k\cos\theta ,
\]
where $h_k$ is well defined on $ST^*K$ since $\cos\theta$ is invariant under $\theta\mapsto-\theta$.  Let
\[
f'_k(\varphi,\psi,\theta)=\bigl(\varphi-k\sin^2\theta,\ \psi+k\theta+k\sin\theta\cos\theta,\ \theta\bigr).
\]

\begin{lemma}\label[lemma]{lem:strict}
The map $f'_k$ is a diffeomorphism of $ST^*K$ with $(f'_k)^*\alpha=\alpha$, and $f'_kf'_\ell=f'_{k+\ell}$.  It is
isotopic to $f_{0,k}$ through diffeomorphisms fixing the base point $(0,0,0)$, so it acts on $H$ by
$s\mapsto s$, $t\mapsto t+ks$, $c\mapsto c$.
\end{lemma}

\begin{proof}
Under $(\psi,\theta)\mapsto(-\psi,-\theta)$ the function $\sin^2\theta$ is invariant and $\sin\theta\cos\theta$ changes sign,
so $f'_k$ is compatible with the relation $(\varphi,\psi,\theta)\sim(\varphi+\pi,-\psi,-\theta)$.  We have
\[
(f'_k)^*\alpha=\cos\theta\,(d\varphi-2k\sin\theta\cos\theta\,d\theta)+\sin\theta\,\bigl(d\psi+k\,d\theta+k(\cos^2\theta-\sin^2\theta)\,d\theta\bigr)=\alpha ,
\]
since the coefficient of $d\theta$ is $k\sin\theta\,(-2\cos^2\theta+1+\cos^2\theta-\sin^2\theta)=0$.  The map $f'_k$ preserves $\theta$
and translates $(\varphi,\psi)$ by a vector depending linearly on $k$, which gives $f'_kf'_\ell=f'_{k+\ell}$.
The maps $(\varphi,\psi,\theta)\mapsto(\varphi-uk\sin^2\theta,\ \psi+k\theta+uk\sin\theta\cos\theta,\ \theta)$, $u\in[0,1]$, form an
isotopy from $f_{0,k}$ to $f'_k$ that fixes the base point, since $\sin\theta=0$ there.
\end{proof}

For $k\neq0$, $f'_k$ does not extend over $DT^*K$, since a diffeomorphism that extends preserves the
kernel $\langle t\rangle$ of $H\to\pi_1(DT^*K)$, while $t+ks$ maps to $r^k\neq1$.

\begin{remark}\label[remark]{rem:nemirovski-maps}
In the flat coordinates above, the two maps of \cite[\S6]{N09} from which the gluing maps are
built are
\[
f_{\mathrm N}(\varphi,\psi,\theta)=(-\varphi,\ \psi+\theta+\pi,\ -\theta-\pi),\qquad
\tau_{\mathrm N}(\varphi,\psi,\theta)=(-\varphi,\ \psi,\ -\theta-\pi),
\]
with $f_{\mathrm N}\circ\tau_{\mathrm N}=f_{0,-1}$ by \cite[Remark~13(ii)]{N09}.  With $h_1=-\cos\theta$ as above, a direct
calculation gives $\tau_{\mathrm N}^*\alpha=\alpha$ and $f_{\mathrm N}^*\bigl(\alpha+\frac12dh_1\bigr)=\alpha+\frac12dh_1$.  Since $dh_1\wedge d\alpha=0$, the forms
$\alpha+\frac u2dh_1$, $u\in[0,1]$, are all contact forms, so by Gray's theorem the contact structures
$\ker\alpha=\xi_{\mathrm{can}}$ and $\ker\bigl(\alpha+\frac12dh_1\bigr)$ are isotopic.  Conjugating $f_{\mathrm N}$ by a diffeomorphism
isotopic to the identity therefore gives a contactomorphism of $\xi_{\mathrm{can}}$ isotopic to $f_{\mathrm N}$.
Composing it with $\tau_{\mathrm N}$ gives a contactomorphism of $\xi_{\mathrm{can}}$ isotopic to $f_{0,-1}$.  Taking powers
and inverses gives the corresponding statement for every Luttinger twist $f_{0,k}$.  So $\tau_{\mathrm N}$ and $f_{\mathrm N}$ are strict
contactomorphisms for two contact forms with isotopic kernels, and \cref{lem:strict} gives
explicit strict representatives for the single contact form $\alpha$.
\end{remark}

\subsection{Rational homology ball fillings}\label{subsec:canonical-fillings}

Let $\lambda_{\mathrm W}=\lambda_{\mathrm{can}}+dF$, where $F$ is a suitable smooth function supported in a small
neighborhood of the zero section, be a standard Weinstein perturbation of the canonical
Liouville structure on $DT^*K$ \cite[Example~11.12(2)]{CE}.  Then
$d\lambda_{\mathrm W}=\omega_{\mathrm{can}}$, so the zero section is Lagrangian, and $\lambda_{\mathrm W}=\lambda_{\mathrm{can}}$ near the boundary.  For a
Legendrian knot $\Lambda$ in $(ST^*K,\xi_{\mathrm{can}})$, let $DT^*K\cup h_\Lambda$ be the Weinstein domain obtained by
attaching a Weinstein $2$--handle along $\Lambda$, and let $(Y_\Lambda,\xi_\Lambda)$ be its boundary.  When we
call a Weinstein domain a Stein filling, we mean a Stein representative of its Weinstein
deformation class \cite{CE,Gompf}.

\begin{theorem}\label[theorem]{thm:qhb-canonical}
Let $\Lambda$ be a Legendrian knot in $(ST^*K,\xi_{\mathrm{can}})$ whose free homotopy class is $(ps+qt)\,c^2$ for
some integers $p$ and $q$ with $q\neq0$, and let $(Y_\Lambda,\xi_\Lambda)$ be the result of Legendrian surgery
along $\Lambda$.  Then, for every $k\in\Z$, the Stein domain obtained by attaching a Weinstein $2$--handle
to $DT^*K$, with its canonical symplectic structure, along $(f'_1)^k(\Lambda)$ is a rational homology ball Stein filling of $(Y_\Lambda,\xi_\Lambda)$ with
fundamental group
\[
\langle r,b\mid brb^{-1}=r^{-1},\ r^{|p+kq|}b^2=1\rangle ,
\]
which is the dicyclic group of order $4|p+kq|$ if $p+kq\neq0$ and the infinite dihedral group if
$p+kq=0$.  In particular,
infinitely many of these fillings are pairwise non-homeomorphic.
\end{theorem}

\begin{proof}
We have $(f'_1)^k=f'_k$ by \cref{lem:strict}, and we write $DT^*K\cup h_{f'_k(\Lambda)}$ for the Stein domain in the
statement.  For an integer $j\ge0$ let
\[
G_j=\langle r,b\mid brb^{-1}=r^{-1},\ r^{j}b^2=1\rangle .
\]
Since $f'_k$ is a contactomorphism of $(ST^*K,\xi_{\mathrm{can}})$, the knot $f'_k(\Lambda)$ is Legendrian, and $f'_k$
maps a standard neighborhood of $\Lambda$ to one of $f'_k(\Lambda)$, so it extends to a contactomorphism
from $(Y_\Lambda,\xi_\Lambda)$ to the result of Legendrian surgery along $f'_k(\Lambda)$.  Hence $DT^*K\cup h_{f'_k(\Lambda)}$,
with the Weinstein structure $\lambda_{\mathrm W}$ on $DT^*K$, is a Weinstein filling, and so a Stein filling,
of $(Y_\Lambda,\xi_\Lambda)$.

By van Kampen's theorem, the fundamental group of $DT^*K\cup h_{f'_k(\Lambda)}$ is the quotient of $\pi_1(K)$ by
the normal closure of the image of the free homotopy class of $f'_k(\Lambda)$.  Since $f'_k$ is
isotopic to $f_{0,k}$ through diffeomorphisms fixing the base point, it maps the conjugacy class of
$g\in H$ to that of $(f'_k)_*(g)$, and by \cref{lem:strict} it fixes $c$ and maps $ps+qt$ to $(p+kq)s+qt$.  So
the free homotopy class of $f'_k(\Lambda)$ is that of $\bigl((p+kq)s+qt\bigr)c^2$, whose image in $\pi_1(K)$ is
$r^{p+kq}b^2$, since $t$ maps to $1$, $s$ to $r$ and $c$ to $b$.  Hence the fundamental group of the
filling is $\langle r,b\mid brb^{-1}=r^{-1},\ r^{p+kq}b^2=1\rangle$.  If $p+kq\ge0$, this is $G_{|p+kq|}$.  If $p+kq<0$, the
substitution $r\mapsto r^{-1}$ preserves the first relation and turns the second one into $r^{|p+kq|}b^2=1$, so
the group is again $G_{|p+kq|}$.  Conjugation by $c$ changes $(p,q)$ to $(-p,-q)$, so the free homotopy
class determines $(p,q)$ only up to sign, but this does not change $|p+kq|$.

It remains to show that the fillings are rational homology balls, and that $G_j$ is the dicyclic group
of order $4j$ if $j>0$ and infinite dihedral if $j=0$.  Let $j\ge0$.  The abelianization of $G_j$ is generated by the
classes $\bar r$ and $\bar b$ with the relations $2\bar r=0$ and $j\bar r+2\bar b=0$.  If $j$ is even, the second
relation becomes $2\bar b=0$, and the abelianization is $\Z/2\oplus\Z/2$.  If $j$ is odd, the second relation
becomes $\bar r=2\bar b$, and then the first one becomes $4\bar b=0$, so the abelianization is $\Z/4$, generated
by $\bar b$.  In both cases the abelianization is finite.  Since the fundamental group of the filling
$DT^*K\cup h_{f'_k(\Lambda)}$ is $G_{|p+kq|}$, its first Betti number is $0$.  The filling is obtained from
$DT^*K$, which deformation retracts onto $K$, by attaching one $2$--handle, so its homology vanishes
in degrees $3$ and $4$, and its Euler characteristic is $\chi(K)+1=1$.  Hence its second Betti number
is $0$ as well, and it is a rational homology ball.

The substitution $r\mapsto r^{-1}$ preserves the relation $brb^{-1}=r^{-1}$ and turns $r^jb^2=1$ into $b^2=r^j$, so
$G_j\cong\langle r,b\mid brb^{-1}=r^{-1},\ b^2=r^j\rangle$.  By \cite[Lemma~5.3]{OzSmooth}, this group is the dicyclic
group of order $4j$ if $j>0$, which for $j=1$ is the cyclic group $\Z/4$, and the infinite dihedral group
if $j=0$.  Hence fillings with different values of $|p+kq|$ have non-isomorphic fundamental
groups.  Since $q\neq0$, the map $k\mapsto p+kq$ is injective, so each value of $|p+kq|$ is attained for at
most two values of $k$, and $|p+kq|\ge|k|\,|q|-|p|$ is unbounded.  Hence $|p+kq|$ takes infinitely many
values, and choosing one $k$ for each value gives infinitely many pairwise non-homeomorphic
fillings.
\end{proof}

\begin{remark}\label[remark]{rem:existence}
Legendrian knots as in \cref{thm:qhb-canonical} exist for every $p\in\Z$ and every nonzero $q\in\Z$, since every free
homotopy class of loops in $ST^*K$ is represented by an embedded knot, and every knot can be
approximated by Legendrian knots.
\end{remark}

\begin{remark}\label[remark]{rem:qhb-compare}
For the two families in \cref{cor:qhb}, the attaching circles have images $r^{-(4mn+1)}b^2$ and
$r^{4mn}b^{-2}$, respectively, in the coordinates of \cite[Lemma~5.2]{OzSmooth}.  The inverse of the
second word is $r^{-4mn}b^2$, since $b^2$ is central, and replacing a relator by its inverse does
not change its normal closure.  So, as in the proof of \cref{thm:qhb-canonical}, the two families have fundamental groups
isomorphic to $G_{4mn+1}$ and $G_{4mn}$, respectively, where $G_j=\langle r,b\mid brb^{-1}=r^{-1},\ r^jb^2=1\rangle$ is the
group appearing in \cref{thm:qhb-canonical}.  Identifying the contact manifolds $(Y_m,\xi_m)$ and $(Y'_m,\xi'_m)$ with
Legendrian surgeries as in \cref{thm:qhb-canonical} would require knowing the classes of $g(d_m)$
and $g(d'_m)$ in the coordinates of \cref{subsec:twists}, which depend on the identification of the
open book with this model.
\end{remark}

\appendix
\section{Proofs of \texorpdfstring{\cref{prop:lantern}}{Proposition 2.3}}\label{app:proofs}

We give two proofs of the relation \eqref{eq:arrangement} below, from which
\cref{prop:lantern} follows at the end of this appendix.  The first explains how
the relation arises, and the second does not depend on it.

\subsection{First proof, via the line arrangement}\label{app:arrangement}

\begin{figure}[ht]
\centering
\begin{tikzpicture}[scale=1.3,line cap=round]
\begin{scope}
\clip (-1.75,-1.75) rectangle (2.35,1.75);
\draw[thick] (-1.6,1.6) -- (1.6,-1.6);              
\draw[thick] (-1.6,1.3) -- (2.3,-0.65);              
\draw[thick] (-1.6,0) -- (2.3,0);                    
\draw[thick] (-1.6,-1.3) -- (2.3,0.65);              
\draw[thick] (-1.6,-1.6) -- (1.6,1.6);              
\foreach \p in {(-1,1),(-1,-1),(0,0),(1/3,-1/3),(1/3,1/3),(1,0)} \fill \p circle (1.4pt);
\end{scope}
\node[right] at (1.6,1.6) {\small $L_1$};
\node[right] at (2.3,0.65) {\small $L_2$};
\node[right] at (2.3,0) {\small $L_3$};
\node[right] at (2.3,-0.65) {\small $L_4$};
\node[right] at (1.6,-1.6) {\small $L_5$};
\node[above] at (-0.88,1.05) {\small $V_2$};
\node[below] at (-0.88,-1.05) {\small $V_1$};
\node at (-0.5,0.2) {\small $V_3$};
\node[below] at (0.36,-0.4) {\small $\beta'$};
\node[above] at (0.36,0.4) {\small $\alpha'$};
\node[above] at (1,0.09) {\small $V_0$};
\end{tikzpicture}
\caption{The five lines of \cref{prop:lantern}, with $y$ scaled by $\frac12$.  Each
intersection point is labeled by the curve it contributes to the relation.}
\label{fig:lines}
\end{figure}
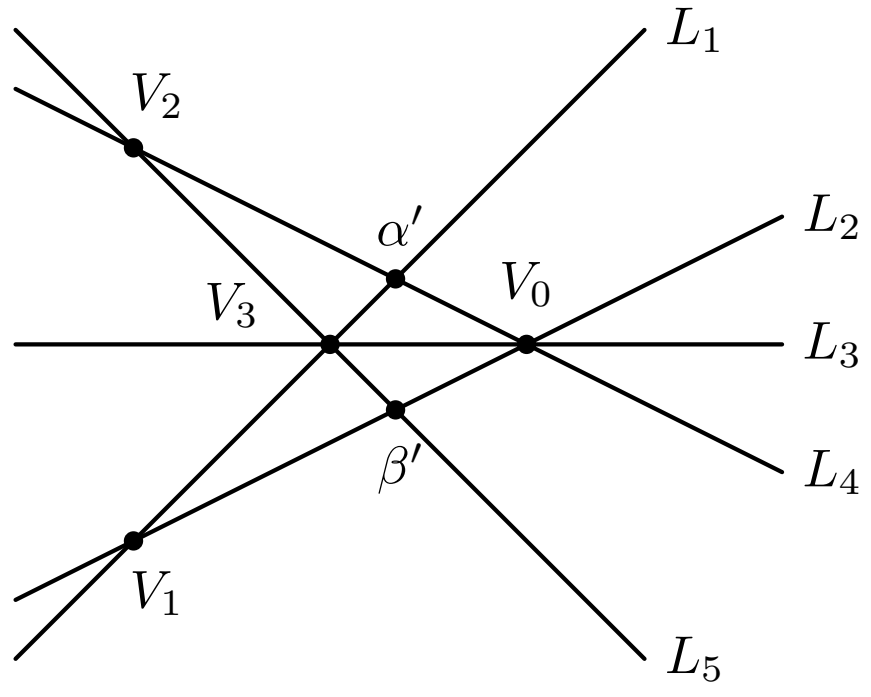

Consider the five lines
\[
L_1\colon y=2x,\quad L_2\colon y=x-1,\quad L_3\colon y=0,\quad L_4\colon y=1-x,\quad
L_5\colon y=-2x
\]
in $\mathbb R^2$, shown in \cref{fig:lines}.  They meet in two triple points,
$p_1=(1,0)\in L_2\cap L_3\cap L_4$ and $p_4=(0,0)\in L_1\cap L_3\cap L_5$, and four
double points, $p_2=L_1\cap L_4$ and $p_3=L_2\cap L_5$ at $x=\frac13$, and
$p_5=L_1\cap L_2$ and $p_6=L_4\cap L_5$ at $x=-1$.  The line $L_3$ contains two
intersection points, and each of the other lines contains three.

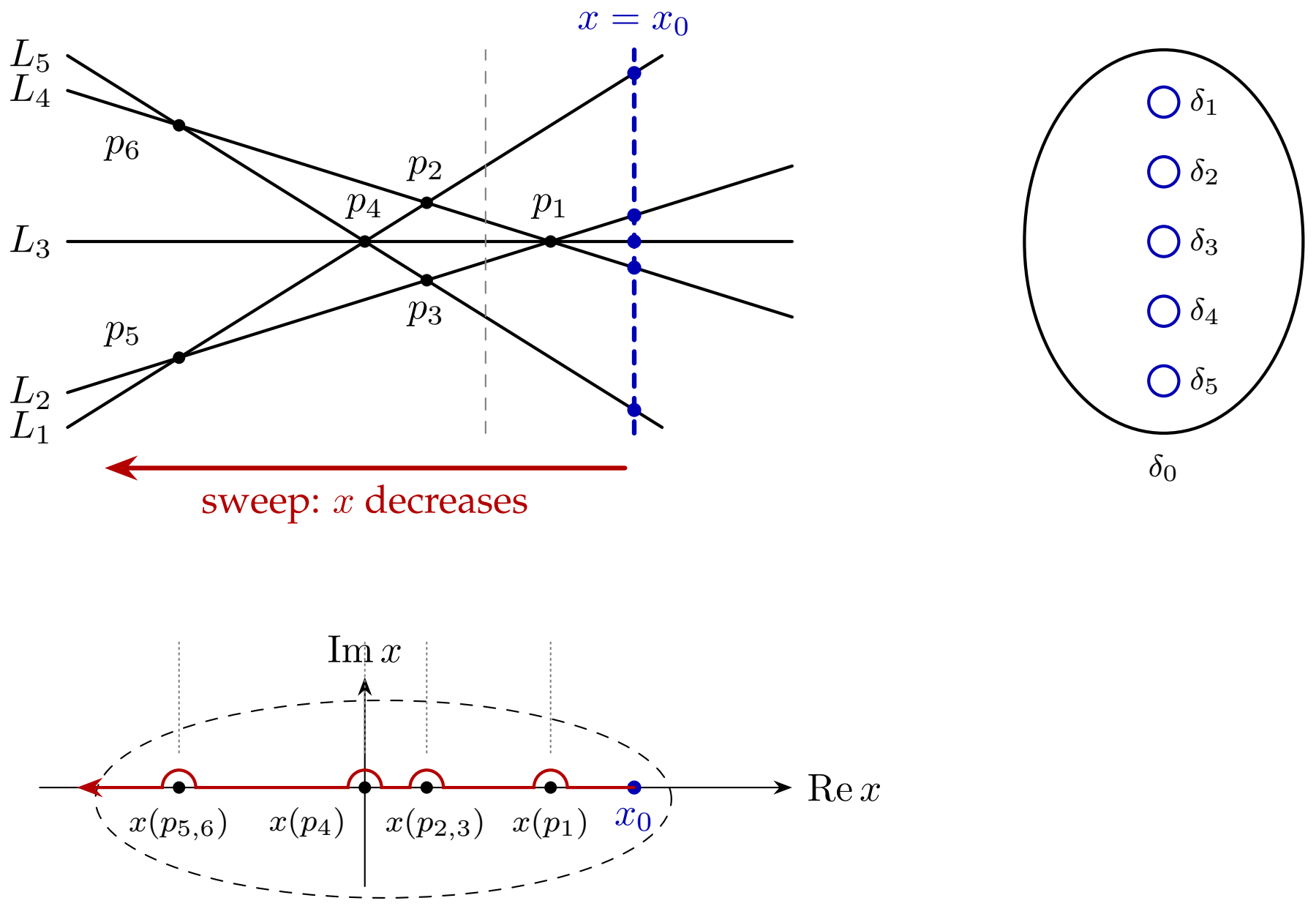
\begin{figure}[ht]
\centering
\begin{tikzpicture}[line cap=round,>=Stealth,scale=1.05]
% ---------- top: the lines and vertical fibers ----------
\begin{scope}[yscale=0.5,xscale=1.6]
\begin{scope}
\clip (-1.75,-3.3) rectangle (2.35,3.3);
\draw[thick] (-1.6,-3.2) -- (1.6,3.2);      % L1: y=2x
\draw[thick] (-1.6,-2.6) -- (2.3,1.3);      % L2: y=x-1
\draw[thick] (-1.6,0) -- (2.3,0);           % L3
\draw[thick] (-1.6,2.6) -- (2.3,-1.3);      % L4: y=1-x
\draw[thick] (-1.6,3.2) -- (1.6,-3.2);      % L5: y=-2x
\end{scope}
\foreach \p in {(1,0),(1/3,2/3),(1/3,-2/3),(0,0),(-1,-2),(-1,2)} \node[circle,fill,inner sep=0pt,minimum size=3.2pt] at \p {};
% reference fiber on the right
\draw[blue!70!black,very thick,dashed] (1.45,-3.3) -- (1.45,3.3);
\foreach \y in {2.9,0.45,0,-0.45,-2.9} \node[circle,fill=blue!70!black,inner sep=0pt,minimum size=3.6pt] at (1.45,\y) {};
\node[blue!70!black,above] at (1.45,3.3) {\small $x=x_0$};
% a second fiber between p_1 and p_2
\draw[gray,dashed] (0.65,-3.3) -- (0.65,3.3);

% sweep arrow
\draw[->,very thick,red!70!black] (1.4,-3.9) -- (-1.4,-3.9) node[midway,below] {\small sweep: $x$ decreases};
\node[left] at (-1.6,-3.2) {\small $L_1$};
\node[left] at (-1.6,-2.6) {\small $L_2$};
\node[left] at (-1.6,0) {\small $L_3$};
\node[left] at (-1.6,2.6) {\small $L_4$};
\node[left] at (-1.6,3.2) {\small $L_5$};
% the fiber over x_0 as a disk with holes
\begin{scope}[shift={(4.3,0)},xscale=1.0,yscale=2]
\draw[thick] (0,0) ellipse (0.75 and 1.65);

\node at (0,-1.95) {\scriptsize $\delta_0$};

\end{scope}
\foreach \y/\l in {2.4/1,1.2/2,0/3,-1.2/4,-2.4/5} {\draw[blue!70!black,thick] (4.3,\y) ellipse (0.081 and 0.26); \node[right] at (4.36,\y) {\scriptsize $\delta_\l$};}
\node at (1,0.6) {\small $p_1$};
\node at (0.33,1.25) {\small $p_2$};
\node at (0.33,-1.25) {\small $p_3$};
\node at (0,0.6) {\small $p_4$};
\node at (-1.3,-1.6) {\small $p_5$};
\node at (-1.3,1.6) {\small $p_6$};
\end{scope}
% the x-plane, directly below: the projection (x,y) -> x
\begin{scope}[shift={(0,-4.7)},xscale=1.6]
\draw[->] (-1.75,0) -- (2.3,0) node[right] {\small $\mathrm{Re}\,x$};
\draw[->] (0,-0.85) -- (0,0.95) node[above] {\small $\mathrm{Im}\,x$};
\draw[dashed] (0.1,-0.1) ellipse (1.55 and 0.85);
\foreach \x in {1,0.333,0,-1} \node[circle,fill,inner sep=0pt,minimum size=3.2pt] at (\x,0) {};
\node[circle,fill=blue!70!black,inner sep=0pt,minimum size=3.6pt] at (1.45,0) {};
\node[blue!70!black,below] at (1.45,-0.05) {\small $x_0$};
\node[below] at (1,-0.05) {\scriptsize $x(p_1)$};
\node[below] at (0.38,-0.05) {\scriptsize $x(p_{2,3})$};
\node[below left] at (-0.02,-0.04) {\scriptsize $x(p_4)$};
\node[below] at (-1,-0.05) {\scriptsize $x(p_{5,6})$};
\draw[red!70!black,thick,->]
  (1.45,0) -- (1.09,0) arc (0:180:0.09 and 0.15) -- (0.423,0) arc (0:180:0.09 and 0.15) -- (0.09,0) arc (0:180:0.09 and 0.15) -- (-0.91,0) arc (0:180:0.09 and 0.15) -- (-1.55,0);
\foreach \x in {1,0.333,0,-1} \draw[gray,densely dotted] (\x,0.3) -- (\x,1.25);
\end{scope}
\end{tikzpicture}
\caption{Top: the sweep in the proof of \cref{prop:lantern}.  The reference
fiber $x=x_0$ (blue) meets $L_1,\ldots,L_5$ in this order from top to bottom, and
it is identified with the disk with holes $\delta_1,\ldots,\delta_5$, shown on the
right.  Moving the vertical fiber from right to left, as indicated by the arrow,
it passes $p_1,\ldots,p_6$.  The gray line is an intermediate fiber.  Bottom: the
$x$--plane, the image of the projection $(x,y)\mapsto x$, with the
$x$--coordinates of the intersection points directly below them, where
$x(p_{2,3})=x(p_2)=x(p_3)$ and $x(p_{5,6})=x(p_5)=x(p_6)$, the path from $x_0$
passing each of them along a small half circle, and a loop around all of them.}
\label{fig:sweep}
\end{figure}

We recall where the relation comes from \cite{Hironaka,BaykurPlanar}.  Regard
the $L_i$ as complex lines in $\mathbb C^2$, and consider the projection
$(x,y)\mapsto x$.  Over a point $x$ of the $x$--plane, the vertical line
$\{x\}\times\mathbb C$ meets the arrangement in five points, one on each $L_i$.
Removing small disks around these points and truncating at large $|y|$, the fiber
becomes the disk with five holes, where the hole $\delta_i$ corresponds to $L_i$.
Over the complement of the $x$--coordinates $x(p)$ of the intersection points,
these fibers form a locally trivial family, so moving $x$ along a loop gives a
mapping class of the fiber.  Near an intersection point $p$, the points of the
lines through $p$ come together, and the monodromy along a small loop around
$x(p)$ is the Dehn twist along a curve enclosing the holes of these lines, up to
twists along the boundaries of these holes, which depend on the framings
\cite[Lemmas~2.1 and~2.2]{Hironaka}.  A
large loop around all the $x(p)$ rotates all five points once around each other,
and its monodromy is a product of boundary twists.  Comparing the two gives the
relation, and capping off all holes but one shows that the exponent of
$t_{\delta_i}$ is $k_i-1$, where $k_i$ is the number of intersection points on $L_i$.
Since the lines are real, all the curves can be compared in one reference fiber
by moving along the real axis of the $x$--plane, passing each $x(p)$ along a small
half circle, as in the bottom of \cref{fig:sweep}.  Along the real axis the five points are real
and ordered from top to bottom, and each half circle reverses the order of the
points of the lines through $p$ by a half twist.

\begin{figure}[htbp]
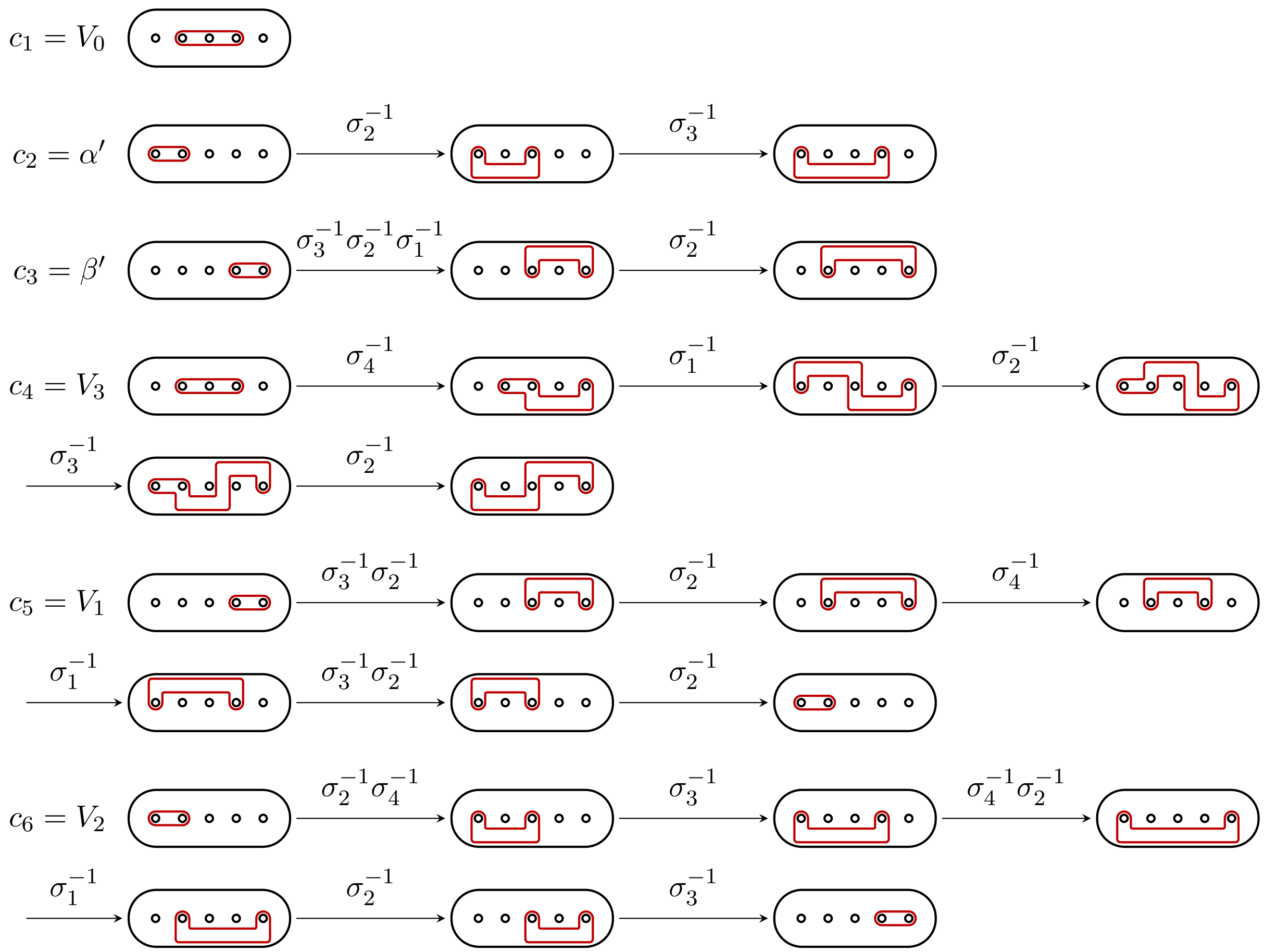

\centering

\caption{The curves $c_1,\ldots,c_6$ of the proof of \cref{prop:lantern}, after
each half twist.}
\label{fig:halftwists}
\end{figure}

We follow the conventions of \cite{BaykurPlanar} for Hironaka's
line-arrangement relation \cite[Theorem~1.2]{Hironaka}.  Choose a vertical fiber
$x=x_0$ to the right of all the intersection points, as in \cref{fig:sweep}.  Its intersections with
$L_1,\ldots,L_5$, ordered from top to bottom, occur in the order $12345$, and we
identify this fiber with the disk with holes $\delta_1,\ldots,\delta_5$, where
$\delta_i$ corresponds to $L_i$.  Move the vertical fiber from right to left, so
that it passes through $p_1,\ldots,p_6$ in this order, where $p_2,p_3$ and $p_5,p_6$
are crossed simultaneously.  These two pairs have the same $x$--coordinate.  The
projection $(x,y)\mapsto x+\varepsilon y$, for small $\varepsilon>0$, separates them,
and since the lines through $p_2$ and $p_3$ are disjoint from each other, as are
those through $p_5$ and $p_6$, the corresponding half twists and Dehn twists commute,
so the order within each pair does not matter.  Upon crossing $p_k$, the lines through $p_k$ form a
consecutive block in the current vertical ordering, and that block is reversed.
The entries in $(i,j)$ below refer to the current positions of the strands, not
to the labels of the lines.
\begin{center}
\begin{tabular}{c|c|c|c|c|c}
point & curve & lines & order before & reversed block & order after\\ \hline
$p_1$ & $V_0$ & $L_2,L_3,L_4$ & $12345$ & $(2,4)$ & $14325$\\
$p_2$ & $\alpha'$ & $L_1,L_4$ & $14325$ & $(1,2)$ & $41325$\\
$p_3$ & $\beta'$ & $L_2,L_5$ & $41325$ & $(4,5)$ & $41352$\\
$p_4$ & $V_3$ & $L_1,L_3,L_5$ & $41352$ & $(2,4)$ & $45312$\\
$p_5$ & $V_1$ & $L_1,L_2$ & $45312$ & $(4,5)$ & $45321$\\
$p_6$ & $V_2$ & $L_4,L_5$ & $45321$ & $(1,2)$ & $54321$
\end{tabular}
\end{center}
The blocks $(1,2)$, $(4,5)$ and $(2,4)$ are reversed by the half twists $\sigma_1$,
$\sigma_4$ and $\tau=\sigma_2\sigma_3\sigma_2$, respectively.  The point $p_k$
contributes the Dehn twist along the curve $c_k$ obtained from the round curve
enclosing the holes in the positions of its block by applying the inverses of the
half twists of $p_1,\ldots,p_{k-1}$, the one of $p_1$ applied last.  Hironaka's
relation \cite[Theorem~1.2]{Hironaka}, in the conventions of \cite{BaykurPlanar},
whose factors occur in the reverse order $p_6,\ldots,p_1$, gives
\[
t_{c_6}t_{c_5}t_{c_4}t_{c_3}t_{c_2}t_{c_1}=t_{\delta_0}\prod_it_{\delta_i}^{\,k_i-1},
\]
where $k_i$ is the number of intersection points on $L_i$.  Here $k_3=2$ and
$k_i=3$ for $i\neq3$, so the right side is $\Delta'$.  By \eqref{eq:artin} and \eqref{eq:artininv},
\begin{align*}
c_1&=x_4x_3x_2=V_0, & c_2&=\tau^{-1}(x_2x_1)=x_4x_1\sim\alpha',\\
c_3&=\tau^{-1}\sigma_1^{-1}(x_5x_4)=x_5x_4x_3x_2x_3^{-1}x_4^{-1}\sim\beta', &
c_4&=\tau^{-1}\sigma_1^{-1}\sigma_4^{-1}(x_4x_3x_2)=x_5x_4x_3x_1x_4^{-1}=V_3,\\
c_5&=\tau^{-1}\sigma_1^{-1}\sigma_4^{-1}\tau^{-1}(x_5x_4)\sim x_2x_1=V_1, &
c_6&=\tau^{-1}\sigma_1^{-1}\sigma_4^{-1}\tau^{-1}\sigma_4^{-1}(x_2x_1)=x_5x_4=V_2.
\end{align*}
The intermediate curves of these computations are shown in
\cref{fig:halftwists}.

Since $V_1$ and $V_2$ are disjoint, as are $\alpha'$ and $\beta'$, the relation becomes
\begin{equation}\label{eq:arrangement}
t_{V_1}t_{V_2}t_{V_3}\,t_{\alpha'}t_{\beta'}\,t_{V_0}
=t_{\delta_0}\,t_{\delta_1}^2t_{\delta_2}^2\,t_{\delta_3}\,t_{\delta_4}^2t_{\delta_5}^2 .
\end{equation}

\Cref{fig:halftwists} shows the curves $c_1,\ldots,c_6$ of the proof of
\cref{prop:lantern}.  The curve $c_1=V_0$ is round, and $c_2,\ldots,c_6$ are obtained
from round curves by applying the half twists one at a time, starting with the
rightmost, using \eqref{eq:artininv}.  Each arrow is labeled
by the half twists applied since the previous picture, and pictures in which the
curve does not change are omitted.  The last picture in each row is the
corresponding curve of \cref{fig:oo-curves}.

\subsection{Second proof, via the action on the fundamental group}\label{app:pi1}

We prove
\eqref{eq:arrangement} directly, for the six curves found above.  This proof uses
neither Hironaka's theorem nor the conventions of the sweep, such as the direction
of the half circles or the framings of the inner boundary components.  Let $L$
denote the left side of \eqref{eq:arrangement}.

First, each of the six Dehn twists is a round twist or the conjugate of a round
twist by an explicit braid.  The curves $V_0$, $V_1$ and $V_2$ are round.  By
\eqref{eq:artin}, $\sigma_1\sigma_4(x_4x_3x_2)=x_4^{-1}x_5x_4x_3x_1$, which is conjugate to
$V_3=x_5x_4x_3x_1x_4^{-1}$, so $V_3=\eta(V_0)$ with $\eta=\sigma_1\sigma_4$.  By
definition $\alpha'=\sigma_1\sigma_2(\beta)$, where $\beta$ is round.  By the computation of
$c_3$ above, $\beta'=\tau^{-1}\sigma_1^{-1}(V_2)=\tau^{-1}(V_2)$, where $\tau=\sigma_2\sigma_3\sigma_2$
is the half twist reversing $\delta_2,\delta_3,\delta_4$, since $\sigma_1^{-1}$ fixes $x_5x_4$.  Hence
\[
t_{V_3}=\eta\,t_{V_0}\,\eta^{-1},\qquad
t_{\alpha'}=\sigma_1\sigma_2\,t_\beta\,(\sigma_1\sigma_2)^{-1},\qquad
t_{\beta'}=\tau^{-1}\,t_{V_2}\,\tau .
\]

Second, we compute the action of $L$ on $\pi_1(\Sigma_0^6,\ast)$.  A right-handed
Dehn twist along a round curve enclosing $\delta_i,\ldots,\delta_j$ conjugates each of
$x_i,\ldots,x_j$ by $x_j\cdots x_i$ and fixes the other generators
(\cref{sec:openbook}), and the half twists act by \eqref{eq:artin} and
\eqref{eq:artininv}.  Applying $t_{V_0}$, $t_{\beta'}$, $t_{\alpha'}$, $t_{V_3}$, $t_{V_2}$ and
$t_{V_1}$ in turn, we find that $L$ conjugates each of $x_1,\ldots,x_5$ by
$x_5x_4x_3x_2x_1$.  The intermediate images are listed at the end of this appendix.  The
curve parallel to $\delta_0$ is the round curve enclosing $\delta_1,\ldots,\delta_5$, so
$t_{\delta_0}$ acts in the same way.  Hence $Lt_{\delta_0}^{-1}$ acts trivially on
$\pi_1(\Sigma_0^6,\ast)$.

Third, a mapping class $f$ of $\Sigma_0^6$ that acts trivially on
$\pi_1(\Sigma_0^6,\ast)$ is a product of twists along $\delta_1,\ldots,\delta_5$.  Indeed,
capping each inner boundary component by a disk with one marked point gives a
homomorphism from the mapping class group of $\Sigma_0^6$ to that of the disk with
five marked points, whose kernel is the free abelian group generated by
$t_{\delta_1},\ldots,t_{\delta_5}$ \cite{FM}.  The image of $f$ acts trivially on the
fundamental group of the disk minus the marked points, which is again the free
group on $x_1,\ldots,x_5$, so it is trivial by Artin's theorem that this action is
faithful \cite{FM}.  Hence $f$ lies in the kernel.  Applied to $f=Lt_{\delta_0}^{-1}$,
this gives
\[
L=t_{\delta_0}\,t_{\delta_1}^{e_1}t_{\delta_2}^{e_2}t_{\delta_3}^{e_3}t_{\delta_4}^{e_4}t_{\delta_5}^{e_5}
\]
for some integers $e_1,\ldots,e_5$.

Finally, we determine the exponents.  Fix $i$, and cap off every inner boundary
component other than $\delta_i$ by a disk.  This gives an annulus with boundary
$\delta_0\cup\delta_i$, whose mapping class group is infinite cyclic, generated by the
twist $t$ along its core.  Under capping, $t_{\delta_0}$ and $t_{\delta_i}$ both map to $t$,
and $t_{\delta_j}$ for $j\neq i$ maps to the identity.  A curve enclosing $\delta_i$
becomes parallel to the core, so its twist maps to $t$, and a curve not enclosing
$\delta_i$ bounds a disk, so its twist maps to the identity.  The six curves enclose
the holes
\[
V_1\colon\{\delta_1,\delta_2\},\quad V_2\colon\{\delta_4,\delta_5\},\quad
V_3\colon\{\delta_1,\delta_3,\delta_5\},\quad \alpha'\colon\{\delta_1,\delta_4\},\quad
\beta'\colon\{\delta_2,\delta_5\},\quad V_0\colon\{\delta_2,\delta_3,\delta_4\},
\]
so $\delta_i$ is enclosed by exactly $k_i$ of them, namely three for $i\neq3$ and two
for $i=3$.  Hence $L$ maps to $t^{k_i}$, and the right side maps to $t^{1+e_i}$.  This
gives $e_i=k_i-1$, that is, $(e_1,\ldots,e_5)=(2,2,1,2,2)$, and \eqref{eq:arrangement}
holds.

We list the action on $\pi_1(\Sigma_0^6,\ast)$ used in the second step.  We apply $t_{V_0}$, $t_{\beta'}$, $t_{\alpha'}$, $t_{V_3}$,
$t_{V_2}$ and $t_{V_1}$ in turn, where $t_{V_3}=\eta\,t_{V_0}\eta^{-1}$ with
$\eta=\sigma_1\sigma_4$, $t_{\alpha'}=\sigma_1\sigma_2\,t_\beta\,(\sigma_1\sigma_2)^{-1}$ and
$t_{\beta'}=\tau^{-1}t_{V_2}\tau$, and each round twist acts as described in
\cref{sec:openbook}.  After each step, the image of $x_j$ is $c_jx_jc_j^{-1}$ for a
word $c_j$, and we list these words, after free reduction.
\begin{list}{}{\setlength{\leftmargin}{1.5em}\setlength{\itemindent}{-1.5em}\setlength{\itemsep}{2pt}\setlength{\topsep}{4pt}}
\item After $t_{V_0}$: $c_{1}=1$, $c_{2}=x_{4}x_{3}$, $c_{3}=c_{4}=x_{4}x_{3}x_{2}$, $c_{5}=1$.
\item After $t_{\beta'}$: $c_{1}=1$, $c_{2}=x_{5}x_{4}x_{3}$, $c_{3}=x_{5}x_{4}x_{3}x_{2}x_{3}^{-1}x_{4}^{-1}x_{5}^{-1}x_{4}$, $c_{4}=x_{5}x_{4}x_{3}x_{2}x_{3}^{-1}x_{4}^{-1}x_{5}^{-1}x_{4}x_{3}$, $c_{5}=x_{5}x_{4}x_{3}x_{2}x_{3}^{-1}x_{4}^{-1}$.
\item After $t_{\alpha'}$: $c_{1}=x_{4}$, $c_{2}=x_{5}x_{4}x_{3}$, $c_{3}=x_{5}x_{4}x_{3}x_{2}x_{3}^{-1}x_{4}^{-1}x_{5}^{-1}x_{4}$, $c_{4}=x_{5}x_{4}x_{3}x_{2}x_{3}^{-1}x_{4}^{-1}x_{5}^{-1}x_{4}x_{3}x_{1}$, $c_{5}=x_{5}x_{4}x_{3}x_{2}x_{3}^{-1}x_{4}^{-1}$.
\item After $t_{V_3}$: $c_{1}=c_{2}=x_{5}x_{4}x_{3}$, $c_{3}=x_{5}x_{4}x_{3}x_{2}x_{1}$, $c_{4}=x_{5}x_{4}x_{3}x_{2}x_{1}x_{4}^{-1}x_{5}^{-1}$, $c_{5}=x_{5}x_{4}x_{3}x_{2}x_{1}x_{4}^{-1}$.
\item After $t_{V_2}$: $c_{1}=c_{2}=x_{5}x_{4}x_{3}$, $c_{3}=c_{4}=c_{5}=x_{5}x_{4}x_{3}x_{2}x_{1}$.
\item After $t_{V_1}$: $c_{1}=x_{5}x_{4}x_{3}x_{2}$, $c_{2}=c_{3}=c_{4}=c_{5}=x_{5}x_{4}x_{3}x_{2}x_{1}$.
\end{list}
In the last line, $c_1=x_5x_4x_3x_2$ may be replaced by $x_5x_4x_3x_2x_1$, since $x_1$
commutes with itself.  So the last line is conjugation by $x_5x_4x_3x_2x_1$, which is
the action of $t_{\delta_0}$.

In either case, since the right side of \eqref{eq:arrangement} is a product of
boundary twists, it is central, and conjugating by $t_{V_0}$ gives
\[
t_{V_0}t_{V_1}t_{V_2}t_{V_3}\,t_{\alpha'}t_{\beta'}
=t_{\delta_0}\,t_{\delta_1}^2t_{\delta_2}^2\,t_{\delta_3}\,t_{\delta_4}^2t_{\delta_5}^2 ,
\]
that is, $B't_{\alpha'}t_{\beta'}=\Delta'$.

\FloatBarrier

\end{document}